%% file: ifacconf.tex
\documentclass{ifacconf}

\usepackage{graphicx}      
\usepackage{natbib}        
\usepackage{booktabs} 
\usepackage[T1]{fontenc}

\usepackage{float}
\usepackage[applemac]{inputenc}
\usepackage{graphics} 
\usepackage{epsfig} 
\usepackage{epstopdf}
\usepackage{amsmath} 
\usepackage{amssymb}  
\usepackage{epstopdf}
\usepackage{graphicx,xcolor}
\usepackage{caption}
\usepackage{subcaption}
\usepackage{courier}
\usepackage{mathrsfs}
\usepackage{pgfplots}
\usepackage{balance}
\usepackage{pgf}
\usepackage{tikz}
\usetikzlibrary{shapes,arrows}
\usetikzlibrary{arrows,automata, decorations.markings}
\DeclareMathOperator*{\argmin}{arg\,min}
\usepackage{tikz}
\usetikzlibrary{shapes,arrows}

\tikzstyle{block} = [draw, fill=blue!10, rectangle, 
    minimum height=4em, minimum width=7em]
\tikzstyle{sum} = [draw, fill=blue!10, circle, node distance=1cm]
\tikzstyle{input} = [coordinate]
\tikzstyle{output} = [coordinate]
\tikzstyle{pinstyle} = [pin edge={to-,thin,black}]

\newtheorem{defi}{Definition}

\DeclareMathOperator*{\minimize}{\textnormal{minimize}}
\begin{document}
\begin{frontmatter}

\title{Compositional Analysis of Hybrid Systems \\ Defined Over Finite Alphabets} 


\author[First]{Murat Cubuktepe} 
\author[First]{Mohamadreza Ahmadi} 
\author[First]{Ufuk Topcu}
\author[Fourth]{Brandon Hencey}

\address[First]{The University of Texas at Austin, Austin, Texas 78712 USA (e-mail: mcubuktepe@utexas.edu, mrahmadi@utexas.edu, utopcu@utexas.edu)}
\address[Fourth]{Air Force Research Laboratory, Dayton, Ohio 45433
USA (e-mail: brandon.hencey@us.af.mil)}

\begin{abstract}
We consider the stability and the input-output analysis problems of a class of large-scale hybrid systems  composed of continuous dynamics coupled with  discrete dynamics defined over finite alphabets, e.g.,  deterministic finite state machines (DFSMs). This class of hybrid systems can be used to model physical systems controlled by software. For such classes of systems, we use a method based on dissipativity theory for compositional analysis that allows us to study stability, passivity and input-output norms. We show that the certificates of the method based on dissipativity theory can be computed by solving a set of semi-definite programs. Nonetheless, the formulation based on semi-definite programs become computationally intractable for relatively large number of discrete and continuous states. We demonstrate that, for systems with large number of states consisting of an interconnection of smaller hybrid systems,  accelerated alternating method of multipliers  can be used to carry out the computations in a scalable and distributed manner.  The proposed methodology is illustrated by an example of a system with 60 continuous states and 18 discrete states.
\end{abstract}

%
%


\end{frontmatter}
\maketitle

\section{Introduction}

Over the past decades, we have witnessed a dramatic increase in research on hybrid and cyber-physical systems~(\cite{AHMADI201437,7867816,7963224,KAMGARPOUR2017177,AHMADI2016118}). Examples of such systems in real-world can be found in robotics~(\cite{7547937}), biological networks~(\cite{Lincoln2004}), and in power systems~(\cite{7726028}).

The literature is rich in analysis and verification methods for hybrid systems ~(\cite{4522632,Alur2011,Livadas1998,4908926}). Despite the available tools for the analysis and verification, scalability  still poses a challenge. Therefore, there has been a surge in compositional analysis techniques. These methods, in general, decompose the analysis problem of a large-scale hybrid system into smaller sub-problems, which can reduce the computational burden significantly. It was shown that dissipativity theory can be used as a tool for decompositional stability and detectability analysis~(\cite{5717501}). This result was further extended in~(\cite{6859327}) to present sufficient conditions for passivity and stability analysis of a class of interconnected hybrid systems with sums of storage functions. Some well-posedness  issues and input-output notions for interconnected hybrid systems were discussed in~(\cite{ricardo52}). Nonetheless, one issue with the above compositional methods for hybrid systems is that they often do not provide a computational framework to find the certificates and rely on ad-hoc analytical techniques to analyze the overall system.


In another vein, several compositional analysis techniques were proposed  based on encoding the hybrid executions in a logic amenable to satisfiability checking (see the survey~\cite{Platzer2011}). \cite{79dd63224} propose a method based on  bounded error approximations
of the hybrid dynamics and the satisfiability checking was carried out using the tool Z3. \cite{7743228} bring forward a verification method at the intersection of software model checking and hybrid systems reachability, which decomposes the discrete and the continuous dynamics.  However, the latter approaches based on SMT formulation are undecidable for general hybrid systems~(\cite{Gao17}) and convergence is not guaranteed.

A large class of discrete systems possess inputs and/or outputs that take values in finite sets. This modeling framework is natural, for example, whenever the actuation takes the form of an on/off or a multi-level switch or when the output is  a binary or a quantized signal. This class of systems is referred to as \textit{systems defined over finite alphabets}. \cite{tarraf2008framework} proposed a robust input-output analysis framework for such systems and \cite{6044564,6882777} applied this framework for controller synthesis.

The contributions of this paper are threefold. First, we propose a stability and an input-output analysis framework for a class of hybrid systems with discrete components defined over a finite alphabet based on Lyapunov and dissipativity analysis. In particular, the discrete dynamics can be in the form of deterministic finite state machines, which can be used to model software. Second, we generalize the methodology proposed in~\cite{dissipativ2016} and \cite{Meissen201555}, wherein a dissipativity-based compositional analysis method was proposed for continuous  systems, to hybrid systems defined over finite alphabets.  Our method decomposes the analysis problem of the overall interconnected hybrid system defined over finite alphabet into smaller local sub-problems for subsystems and  takes advantage of a global storage function, which is the sum of local storage functions. Thirdly, to carry out the computation in a distributed manner, we use accelerated ADMM (\cite{goldstein2014fast}), which is a variant of ADMM (\cite{boyd2011distributed}). To use accelerated ADMM, which has a faster convergence rate compared to ADMM, we utilize \emph{smoothing techniques} (\cite{nesterov2005smooth,becker2011templates}), which has been used to improve the convergence rate of similar first order methods. We also discuss the effects of \emph{restarting} accelerated ADMM, which has shown to improve the convergence rate of similar accelerated algorithms (\cite{nesterov2013gradient,becker2011templates}). We illustrate the proposed method by a numerical example.

This paper is structured as follows. In the following section, we present the problem formulation. In Section 3, we define notions of hybrid Lyapunov and storage functions for hybrid systems defined over finite alphabets. In Section~\ref{sec:interc}, we propose a method based on dissipativity for compositional analysis of the class of hybrid systems under study. The proposed methodology is illustrated by an example in Section~\ref{sec:examples}. Finally, Section~\ref{sec:conclude} concludes the paper and gives directions for future research. 

\textbf{Notation:}
$\mathbb{R}_{\ge 0}$ denotes the set $[0,\infty)$. $\| \cdot \|$ denotes the Euclidean vector norm on $\mathbb{R}^n$. The set of integers are denoted by $\mathbb{Z}$.  For a function $f:A\to B$, $f \in L_p(A, B)$, $1\le p  < \infty$, implies that $\left( \int_A |f(t)|^p dt \right)^{\frac{1}{p}} < \infty$ and $\sup_{t \in A} |f(t)| < \infty$ for $p=\infty$. Equivalently, for a discrete signal $s: \mathcal{A} \to \mathcal{B}$, $s \in l_p(\mathcal{A}, \mathcal{B})$, $1\le p  < \infty$, implies that $\left( \sum_{\mathcal{A}} | s(n)|^p  \right)^{\frac{1}{p}} < \infty$ and $\sup_{n \in A} |s(n)| < \infty$ for $p=\infty$.  For symmetric matrices $A_1,\ldots,A_n$, $diag(A_1,\ldots,A_n)$ denotes the diagonalized matrix
$$
diag(A_1,\ldots,A_n) = \begin{bmatrix} A_1 & 0 & 0 \\ 0 & \ddots & 0 \\ 0 & 0 & A_n \end{bmatrix}.
$$
For a vector $s \in \mathbb{R}^{n_s}$, $s \equiv 0$ denotes the element-wise equality to zero. If $x \in X$, then $f \in \mathcal{C}^1(X)$ implies that the function $f$ is continuously differentiable in $x$.

\section{Problem Formulation}



Formally, we consider the following class of hybrid systems
\begin{eqnarray} \label{eq:HS}
\mathcal{G}: \begin{cases}
\mathcal{C}: \begin{cases} \dot{x}(t) = f\left(x(t),w(t);p(t)\right) \\ y(t) = h\left(x(t);p(t)\right)  \end{cases}\\
\mathcal{D}: \begin{cases} q(t^{+}) = g\left(q(t),u(t),x(t)\right) \\ p(t) = l\left(q(t),u(t)\right)  \end{cases}\\
q(t_0) = q_0,~x(t_0)=x_0.
\end{cases}
\end{eqnarray}
where $x \in \mathbb{R}^n$ and $q \in \mathscr{Q} \subset \mathbb{N}$ represent continuous and discrete states. In the continuous module $\mathcal{C}$, $f(\cdot,\cdot;p):\mathbb{R}^n \times \mathbb{R}^m \to \mathbb{R}^n$, $f(0,w;p) \equiv 0$, $\forall (w,p) \in \mathscr{W} \times \mathscr{P}$, is a family of mappings with index $p \in \mathscr{P}\subset \mathbb{Z}$ and similarly $h(\cdot;p):\mathbb{R}^n  \to \mathbb{R}^{n_y}$ is a family of output mappings. $y \in \mathbb{R}^{n_y}$ and $w \in \mathbb{R}^m$ are the continuous outputs and inputs, respectively. In the discrete module $\mathcal{D}$, $g:\mathscr{Q} \times \mathscr{U} \times \mathbb{R}^n  \to \mathscr{Q}$ and $l:\mathscr{Q} \times \mathscr{U} \to \mathscr{Q}$. $p \in \mathscr{P}$ and $u \in \mathscr{U}\subset \mathbb{Z}^{n_u}$ are the discrete outputs and inputs, respectively. The sets associated with the discrete module $\mathscr{Q}$, $\mathscr{P}$ and $\mathscr{U}$ are assumed to be finite. In the sequel, we abuse the notation and use $q^{+}$ to represent $q(t^{+})$.

The discrete module $\mathcal{D}$ can characterize a rich class of systems defined over finite alphabets~(\cite{tarraf2008framework}) and can be used to model systems ranging from quantizers to deterministic finite state machines. The hybrid system~\eqref{eq:HS} can also be studied in the context of hybrid automata~(\cite{Henzinger1996}). However, note that for hybrid automata, analysis tools such as Lyapunov functions or storage functions are not available in general.
  
  \section{Stability and Dissipativity Analysis}
  
  We can study the input-output and stability properties of system~\eqref{eq:HS} by using a dissipativity-type and Lyapunov-type argument, respectively. To this end, we use the notion of hybrid Lyapunov or storage function, which is described as follows.
  
  \begin{defi} [Hybrid Lyapunov Function]
  A function $V: \mathbb{R}^n \times \mathscr{Q} \to \mathbb{R}_{\ge 0}$ such that $V(0,q)=0,~\forall q \in \mathscr{Q}$, and $V \in \mathcal{C}^1(\mathbb{R}^n)$ is called a \emph{hybrid Lyapunov function} for  system~\eqref{eq:HS} with $u \equiv 0$ and $w \equiv 0$, if it satisfies the following inequalities
  \begin{equation} \label{inq1stab}
 V(x,q)>0,\quad \forall x \in \mathbb{R}^n\setminus \{0\},~\forall q \in \mathscr{Q},
  \end{equation}
    \begin{equation}\label{inq2stab}
  \left( \frac{\partial V(x,q)}{\partial x} \right)^T f(x,0;p) < 0,\quad \forall x \in \mathbb{R}^n,~~\forall q \in \mathscr{Q},~~\forall p \in \mathscr{P},
  \end{equation}
  and
      \begin{equation}\label{popos}
  {V}(x,q^{+})- V(x,q) \le 0,\quad \forall x \in \mathbb{R}^n, \forall q \in \mathscr{Q}.
  \end{equation}
  \end{defi}
  
 \begin{thm}
 The hybrid system~\eqref{eq:HS} is asymptotically stable, i.e., $\lim_{t \to \infty} x(t) = 0,~~\forall q \in \mathscr{Q}$, if there exists a hybrid Lyapunov function.
 \end{thm}
 
 \begin{pf}
See Appendix~\ref{app1}.
 \end{pf}
 
   \begin{defi} [Hybrid Storage Function]
  A function $V: \mathbb{R}^n \times \mathscr{Q} \to \mathbb{R}_{\ge 0}$ such that $V \in \mathcal{C}^1(\mathbb{R}^n)$ is called a \emph{hybrid storage function} for  system~\eqref{eq:HS}, if it satisfies the following inequalities 
  \begin{equation}
  V(x,q) \ge 0,\quad \forall x \in \mathbb{R}^n, \forall q \in \mathscr{Q},\label{eq:HS_storage1}
  \end{equation}
    \begin{multline}
  \left( \frac{\partial V(x,q)}{\partial x} \right)^T f(x,w;p)\le W_c(w,y),\\ \forall x \in \mathbb{R}^n,~~\forall q \in \mathscr{Q},~~\forall p \in \mathscr{P}\label{eq:HS_storage2}
  \end{multline}
  and
      \begin{equation}
  {V}(x,q^{+})- V(x,q) \le W_d(u,p),\quad \forall x \in \mathbb{R}^n, \forall q \in \mathscr{Q}\label{eq:HS_storage3}
  \end{equation}
  where the integrable functions $W_c: \mathbb{R}^m \times \mathbb{R}^{n_y} \to \mathbb{R}$ and $W_d:\mathscr{Q} \times \mathscr{U} \to \mathbb{R}$ are the continuous and the discrete supply rates, respectively.
  \end{defi}
  
   \begin{thm}
 The hybrid system~\eqref{eq:HS} is dissipative with respect to the supply rates $W_c$ and $W_d$,  if there exists a hybrid storage function.
 \end{thm}
 
   \begin{pf}
 The dissipativity of the continuous dynamics is standard and follows from integrating~\eqref{eq:HS_storage2}. The dissipativity of the discrete dynamics follows from Theorem 3 and Lemma 2 in~\cite{tarraf2008framework}.
    \end{pf}
 
\begin{figure}[H]
\centering
\includegraphics[width=6cm]{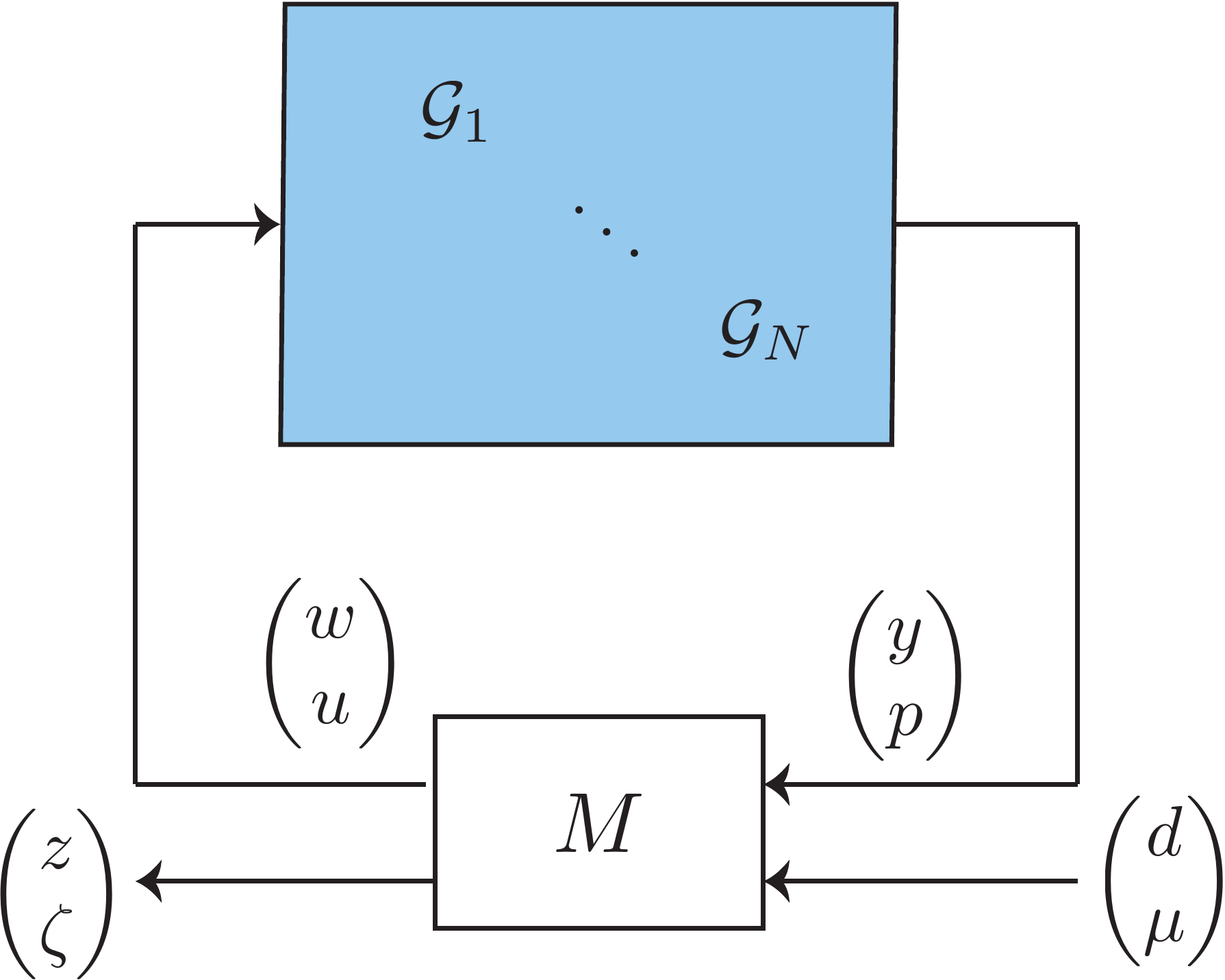}
\caption{Interconnected system with hybrid inputs $(d,\mu)^T$ and hybrid outputs $(z,\zeta)^T$.}\label{fig1}
\end{figure}
  
  \section{Interconnection of Hybrid Systems Defined Over Finite Alphabets}\label{sec:interc}
  We consider interconnected systems as illustrated in Fig.~\ref{fig1}, where the subsystems $\{ \mathcal{G}_i \}_{i=1}^N$ are known and have dynamics in the form of~\eqref{eq:HS}. We associate each subsystem with a set of functions $\{ f_i,h_i,g_i,l_i\}$ and $x_i \in  \mathbb{R}^{n_i}$, $q_i \in \mathscr{Q}_i$, $w_i \in  \mathbb{R}^{n_w^i}$, $u_i \in \mathscr{U}_i$, $y_i \in \mathbb{R}^{n_y^i}$ and $p_i \in \mathscr{P}_i$. The static interconnection is characterized by a matrix $M $ where $n = \sum_{i=1}^N n_i$, $n_w = \sum_{i=1}^N n_w^i$ and $n_y = \sum_{i=1}^N n_y^i$. That is, $M$ satisfies 
\begin{equation} \label{eq:interconnection}
\begin{bmatrix} w \\ z \\ u \\ \zeta \end{bmatrix} = M \begin{bmatrix} y \\ d \\ p \\ \mu \end{bmatrix},
\end{equation}
where $d \in \mathbb{R}^{n_d}$ and $z \in \mathbb{R}^{n_z}$ are the continuous exogenous inputs and outputs, respectively. Similarly, $\mu \in \mathscr{M} \subset \mathbb{Z}^{n_{\mu}}$ and $\zeta \in \mathscr{Z} \subset \mathbb{Z}^{n_{\zeta}}$ are the discrete exogenous inputs and outputs, respectively. We assume this interconnection is well-posed, i.e., for all $d \in L_{2e}$ and initial condition $x(0) \in \mathbb{R}^n$, there exist unique $z,w,y \in L_{2e}$  that causally depend on $d$ for all $p \in \mathscr{P} \subset \mathbb{Z}^{n_p}$. Furthermore, we define
$$
M=
\left[
\begin{array}{c|c}
M_{c}&  \\
\hline
& M_{d}
\end{array}
\right],
$$
where $M_{c} \in \mathbb{R}^{n_w + n_y} \times \mathbb{R}^{n_w + n_y}$ and $M_{d} \in \mathbb{Z}^{n_\mu + n_\zeta} \times \mathbb{Z}^{n_\mu + n_\zeta}$. 
The continuous local  and global supply rates, $W_c^i(w_i,y_i)$ and $W_c(d,z)$, respectively, are defined by quadratic functions. That is, 
\begin{equation} \label{fsfFF}
W_c(d,z) = \begin{bmatrix}  d \\ z \end{bmatrix}^T S \begin{bmatrix}  d \\ z \end{bmatrix}.
\end{equation}
where $S \in \mathbb{R}^{n_d+n_z} \times \mathbb{R}^{n_d+n_z}$.

Analogously, for discrete subsystems, we have $W_d^i(p_i,u_i)$ and $W_d(p,u)$ given by quadratic functions
\begin{equation} \label{fsfFFd}
W_d(p,u) = \begin{bmatrix}  u \\ p \end{bmatrix}^T R \begin{bmatrix}  u \\ p \end{bmatrix}.
\end{equation}
where $R \in \mathbb{Z}^{n_u+n_p} \times \mathbb{Z}^{n_u+n_p}$ is a symmetric matrix.

We next show that certifying the dissipativity of an overall interconnected system can be concluded, if each of the subsystems satisfy the local dissipativity property. Let 
\begin{multline}
\mathcal{L}_i =  \bigg \{ (S_i,~R_i) \mid   \text{$\mathcal{G}_i$ is dissipative w.r.t.}\\\text{$\begin{bmatrix} w_i \\ y_i \end{bmatrix}^T S_i \begin{bmatrix} w_i \\ y_i \end{bmatrix}$},~ \text{and}~ \begin{bmatrix}  u_i \\ p_i \end{bmatrix}^T R_i \begin{bmatrix}  u_i \\ p_i \end{bmatrix} \bigg\},
\end{multline}
\begin{equation} \label{sddccxcxglobal}
\mathcal{L}_c = \bigg\{ S, \{S_i\}_{i=1}^N \mid \begin{bmatrix} M_c \\ 1_{n_y}  \end{bmatrix}^T P_c^T Q_c P_c \begin{bmatrix} M_c \\ 1_{n_y}   \end{bmatrix} < 0 \bigg\},
\end{equation}
and
\begin{equation} \label{sddccxcxglobald}
\mathcal{L}_d = \bigg\{ R, \{R_i\}_{i=1}^N \mid \begin{bmatrix} M_d \\ 1_{n_p}  \end{bmatrix}^T P_d^T Q_d P_d \begin{bmatrix} M_d \\ 1_{n_p}   \end{bmatrix} < 0 \bigg\},
\end{equation}
wherein $Q_c = diag(S_1,\ldots,S_N,-S)$, \\$Q_d = diag(R_1,\ldots,R_N,-R)$, and $P_c$ and $P_d$ are permutation matrices defined by
\begin{equation}
\begin{bmatrix} w_1 \\ y_1 \\ \vdots \\ w_N \\ y_N \\d \\z \end{bmatrix} = P_c \begin{bmatrix} w \\ z  \\y \\d \end{bmatrix},~\text{and}~\begin{bmatrix} u_1 \\ p_1 \\ \vdots \\ u_N \\ p_N \\ \mu \\ \zeta \end{bmatrix} = P_d \begin{bmatrix} u \\ \zeta  \\ p \\ \mu \end{bmatrix} .
\end{equation}

\begin{prop}\label{prop1}
Consider the interconnection of $N$ subsystems as given in~\eqref{eq:interconnection} with the global supply rates~\eqref{fsfFF} and \eqref{fsfFFd}. If there exists $\{ S_i \}_{i=1}^N$ and $\{ R_i \}_{i=1}^N$ satisfying
\begin{equation}\label{ddsdsdsdsdsccdcdxx}
(S_i,R_i) \in \mathcal{L}_i,~~i=1,\ldots,N,
\end{equation}
and
\begin{equation} \label{eq:global}
(S_1,\ldots, S_N, -S) \in \mathcal{L}_c,
\end{equation}
\begin{equation} \label{eq:globald}
(R_1,\ldots, R_N, -R) \in \mathcal{L}_d,
\end{equation}
then the interconnected system is dissipative with respect to the global supply rates $W_d$ and $W_c$. A storage function certifying  global dissipativity is $V(x,q) = \sum_{i=1}^N V_i(x_i,q_i)$, where $V_i$ is the storage function certifying dissipativity of subsystem $i$ as in $\mathcal{L}_i$.
\end{prop}
\begin{pf}
See Appendix~\ref{app3}.
\end{pf}

Proposition 3 provides the means to decompose the analysis of interconnected subsystems to smaller problems that are computationally more amenable. However, even the above discussed decompositional analysis method can be computationally involved for large-scale hybrid systems. In the next section, we propose a method based on accelerated ADMM to carry out such computations in a distributed manner. 

 It is worth mentioning \cite{Oehlerking2009} proposed a compositional method based on graph-based reasoning rather than dissipativity for hybrid  systems that are not subject to inputs and outputs and are described by a hybrid automaton. Whereas, our formulation allows for (both discrete and continuous) inputs and outputs and is pertained to hybrid systems composed of continuous dynamics and a discrete subsystem defined over finite alphabets.  In addition,  \cite{Oehlerking2009} did not bring forward a method based on distributed optimization  to address the computations in the case of large-scale systems. Such distributed computational method is discussed in the next section.

\section{Computational Formulation Using Accelerated ADMM}\label{sec:admm}
For small-scale systems, we can solve the optimization problem outlined in Proposition 3 using publicly available SDP solvers like MOSEK (\cite{andersen2012mosek}), SeDuMi (\cite{sturm1999using}) or SDPT3 (\cite{toh1999sdpt3}). But, these SDP solvers do not scale well for larger problems, as they use interior point methods, which requires solving a system of equations in each iteration. However, the structure in our problem allows us to decompose the constraints in \eqref{eq:HS_storage1}, \eqref{eq:HS_storage2} and \eqref{eq:HS_storage3}, leading to a distributed algorithm. Specifically, the ADMM (\cite{boyd2011distributed}) approach decomposes the convex optimization problems into a set of smaller problems. A generic convex optimization problem
\begin{align}
 \text{minimize} &\quad  F(b)\nonumber\\
 \text{subject to} &\quad  b\in \mathit{C},
\end{align}
where $b \in \mathbb{R}^n$, $F$ is a convex function, and $C$ is a convex set, can be written in ADMM form as 

\begin{align}
 \text{minimize} &\quad  F(b)+G(v)\nonumber \\ 
 \text{subject to} &\quad  b=v, \label{eq:admm2}
\end{align}
where $G$ is the indicator function of $\mathit{C}$.

The problem we want to find a compositional formulation can be given as

\begin{eqnarray}
&\displaystyle \minimize_{\{S_i\}_{i=1}^N,\{R_i\}_{i=1}^N, \{V_i\}_{i=1}^N} \eta & \nonumber \\
&\text{subject to \eqref{eq:HS_storage1},\eqref{eq:HS_storage2},\eqref{eq:HS_storage3}, \eqref{eq:global}, and \eqref{eq:globald}}& \label{eq:globalobj}
\end{eqnarray}
which is outlined in Proposition 3.

For example, $\eta$ can be the upper-bound on the continuous induced norm  $\frac{\|z\|_{L_2}}{\|d\|_{L_2}}$ or the discrete induced norm $\frac{\|\zeta\|_{l_2}}{\|\mu\|_{l_2}}$ that we wish to minimize. Using the above form, the  problem in \eqref{eq:globalobj} can be written in ADMM form with $F(b)$ is defined as sum of $\eta$ and the indicator function of \eqref{eq:HS_storage1},\eqref{eq:HS_storage2} and \eqref{eq:HS_storage3}, and $G(v)$ is defined as the indicator function of \eqref{eq:global} and \eqref{eq:globald}. Then, the scaled form of ADMM algorithm for problem in \eqref{eq:admm2} is

\begin{align*}
& b^{k+1} = \argmin_{b} F(b) + (\rho /2)\vert\vert b-v^k+s^k\vert\vert^2_2,\\
& v^{k+1}=\argmin_{v} G(v) + (\rho /2)\vert\vert b^{k+1}-v+s^k\vert\vert^2_2,\\
& s^{k+1}=s^k+b^{k+1}-v^{k+1},
\end{align*}
where $b$ and $v$ are the vectorized form of the matrices $\{S_i\}_{i=1}^N$, $\{V_i\}_{i=1}^N$, $\{R_i\}_{i=1}^N$, $z$ is the scaled dual variable and $\rho > 0$ is the \emph{penalty parameter}, which is a penalty for primal infeasibility, i.e penalty of not satisfying the constraint $b=v$. As $F(b)$ is separable for each subsystem, the ADMM algorithm can be parallelized as follows:

\begin{align*}
& b_i^{k+1} = \argmin_{b_i} F_i(b) + (\rho /2)\vert\vert b_i-v_i^k+s_i^k\vert\vert^2_2,\\
& v^{k+1}=\argmin_{v} G(v) + (\rho /2)\vert\vert b^{k+1}-v+s^k\vert\vert^2_2,\\
& s^{k+1}=s^k+b^{k+1}-v^{k+1},
\end{align*}

Under mild assumptions, the ADMM algorithm converges \cite{boyd2011distributed}, but the convergence is only asymptotic in general, therefore it may require many iterations to achieve sufficient accuracy.

\subsection{Accelerated ADMM}

Several algorithms  (\cite{hale2008fixed, becker2011nesta, nesterov2013gradient, beck2009fast, chen2012fast}) shows that acceleration schemes can improve the performance significantly. These methods achieve $O(\frac{1}{k^2})$ convergence after $k$ iterations, which is shown to be optimal for a first order method \cite{nesterov1983method}. However, they usually require the function $F(b)$ to be differentiable with a known Lipschitz constant on the $\nabla F(b)$, which does not exist when the problem is constrained. For the case when $F(b)$ or $G(v)$ is not strongly convex or smooth, smoothing approaches have been used \cite{nesterov2005smooth,becker2011templates} to improve convergence. However, to the best of our knowledge, these methods have not been applied in compositional analysis.

Consider the following perturbation of the problem in \eqref{eq:globalobj}

\begin{eqnarray}
&\displaystyle \minimize_{\{S_i\}_{i=1}^N,\{R_i\}_{i=1}^N, \{V_i\}_{i=1}^N} \eta +\ell ~D_i(S_i,R_i,V_i,)& \nonumber \\
&\text{subject to \eqref{eq:HS_storage1},\eqref{eq:HS_storage2},\eqref{eq:HS_storage3}, \eqref{eq:global}, and \eqref{eq:globald}}& \label{eq:globalobj_approx}
\end{eqnarray}
for some fixed smoothing parameter $\ell >0$ and a strongly convex function $D$ that satisfies

\begin{equation}
D(b)\geq D(b_0) + \dfrac{1}{2}\vert\vert b-b_0\vert\vert^2_2
\end{equation}
for some point $b_0$. Specifically, we choose $D_i=\lVert S_i\rVert_F+\lVert V_i\rVert_F+\lVert R_i\rVert_F$, where  $\| \cdot \|_F$ is the Frobenius norm.  For some problems, it is shown that for small enough $\ell$, the approximate problem \eqref{eq:globalobj_approx} is equivalent to the original problem \cite{becker2011templates}. 

When $F(b)$ and $G(v)$ are strongly convex, the ADMM algorithm can be modified with an acceleration step to achieve $O(\frac{1}{k^2})$ convergence after $k$ iterations  \cite{goldstein2014fast}. Then, the accelerated ADMM algorithm is

\begin{align*}
& b_i^{k} = \argmin_{b_i} F_i(b) + (\rho /2)\vert\vert b_i-\bar{v}_i^k+\bar{s}_i^k\vert\vert^2_2,\\
& v^{k}=\argmin_{v} G(v) + (\rho /2)\vert\vert b^{k}-v+\bar{s}^k\vert\vert^2_2,\\
& s^{k}=\bar{s}^k+b^{k}-v^{k},\\
& \alpha_{k+1}= \dfrac{1+\sqrt{1+4\alpha^2_k}}{2}\\
& \bar{v}^{k+1}=v^k+\dfrac{\alpha_k -1}{\alpha_{k+1}}(v^k-v^{k-1})\\
& \bar{s}^{k+1}=s^{k}+\dfrac{\alpha_k -1}{\alpha_{k+1}}(s^k-s^{k-1}),
\end{align*}
where $\rho$ is a positive constant that satisfies $\rho\leq \ell$ to make sure the accelerated ADMM converges with $O(\dfrac{1}{k^2})$ rate, and~$\alpha_1=1$.

Note that $\alpha$ update can be carried out in parallel while achieving $O(\frac{1}{k^2})$ convergence, which cannot be achieved by the standard ADMM or accelerated proximal methods if there are constraints in the problem.

In general, we do not have access to the Lipschitz constant or strongly convexity parameter in the feasible region of the subproblems because of the constraints, which may reduce the performance of the accelerated method \cite{becker2011templates}. One approach to deal with the case of unknown Lipschitz constant or strongly convexity parameter is so-called \emph{restart} method, which is used in \cite{ nesterov2013gradient,becker2011templates}, and it is shown that restart methods can improve the convergence rate significantly. To apply the method, we \emph{restart} the algorithm, i.e, we set the acceleration parameter $\alpha_k=1$ after a certain number of iterations while using the point in iteration $k$ as the starting point for the restart, which resets the acceleration parameter, and reruns the accelerated ADMM algorithm from the next starting point. Examples in \cite{becker2011templates} show that the restarting methods can greatly improve performance, but they note that the restart method requires tuning for different problems to optimize the performance, i.e, the performance can vary significantly with different restart schemes.

\section{Numerical Experiments}\label{sec:examples}

In this section, we illustrate the proposed distributed analysis method with a large scale example, where we compare the convergence rate of ADMM with accelerated ADMM and several restart methods. We implemented both standard ADMM and accelerated ADMM algorithms in MATLAB using the CVX toolbox (\cite{grant2008cvx}) and MOSEK (\cite{andersen2012mosek}) to solve SDP problems.

\subsection{Example}
In order to test and compare the different forms of ADMM methods discussed in Section~\ref{sec:admm}, we randomly generated $N=6$ subsystems with linear continuous dynamics, each with 10 continuous states, 2 continuous inputs and outputs, and with 3 discrete modules. The continuous dynamics of each subsystem is characterized by the equations \eqref{eq:HS}, with the maximum real part for an eigenvalue of $A_i$ is normalized to~$-2$, and all of the subsystems are connected to 3 other subsystems, which forms the interconnection matrix $M$. In total, the system we consider has 60 continuous states with 3 different dynamics depending on the discrete module, which makes the centralized approaches impractical as the number of semidefinite variables and constraints grow large. We note that, MOSEK run into numerical problems while solving the analysis problem with a centralized approach, therefore we only include comparisons between different compositional methods.

\begin{figure}[H]
\centering
\begin{tikzpicture}[->,>=stealth',shorten >=2pt,auto,node distance=4cm,
                    semithick]
  \tikzstyle{every state}=[text=black]

  \node[state] (A)                    {$q_3$};
  \node[state]         (B) [above right of=A] {$q_1$};
  \node[state]         (C) [right of=A] {$q_2$};

  \path (A) edge  [bend left]   node {$u=0,1$} (B)
        (B) edge [loop above] node {$u=0$} (B)
         (B)   edge   [bend left]            node {$u=1$} (C)
        (C) edge    [loop right]           node {$u=0$} (C)
         (C)   edge  node {$u=1$} (A);
\end{tikzpicture}

\caption{$\mathcal{D}_i$ in the Example.}
\label{fig:FSM}

\end{figure}
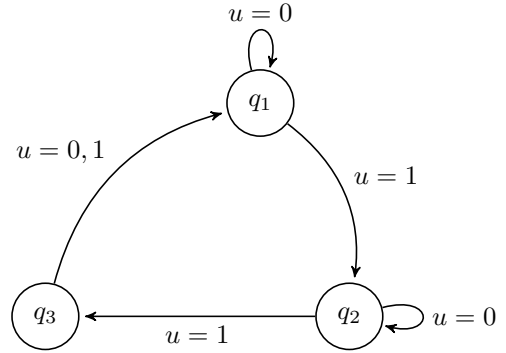

The discrete module $\mathcal{D}_i$ is shown in Figure~\ref{fig:FSM} with three states $q_1, q_2,$ and $q_3$, inputs $\mathscr{U}=\lbrace 0, 1 \rbrace$ and outputs $\mathscr{P}=\lbrace 0, 1 ,2 \rbrace.$ The output function $p$ is defined as

\begin{align*}
p(q_i,u): \begin{cases} u, & \text{for} ~ i =1,\\
1-u, & \text{for} ~ i =2,\\
0 , & \text{for} ~ i =3.
\end{cases}
\end{align*} 


We apply the compositional approach underlined by the problem in \eqref{eq:globalobj_approx} to find the minimum  induced $L_2$-norm between the output and input of randomly chosen systems.

Since the continuous dynamics are linear, we consider quadratic Lyapunov functions for subsystems. For each subsystem, let

\begin{align*}
V(x_i)=\begin{pmatrix}
x_i\\
q_i
\end{pmatrix}^T \begin{bmatrix} P_i & r_i\\ r_i^T & \lambda_i \end{bmatrix} \begin{pmatrix}
x_i\\
q_i
\end{pmatrix}.
\end{align*}

The iterative methods were initialized using $V^0_i=S^0_i=R_i^0=U_i^0=I$. For each method, we plot the norm of \emph{primal residual} in Figure~\ref{primal}, which is defined as $r^{k}=b^k-v^k$, and it is the residual for primal feasibility. Also, we show the  norm of the \emph{dual residual} $s^k=\rho(v^k-v^{k-1})$ in Figure~\ref{dual}, which can be viewed as a residual for the dual feasibility condition. 

\input{ADHS_hybrid_primal.tex}

\input{ADHS_hybrid_dual.tex}

We can see that accelerated ADMM achieves superior convergence in primal and dual residuals compared to ADMM. In Fig.~\ref{primal}, we can see that ADMM can only achieve accuracy up to $10^{-1}$ after 200 iterations, which is not a sufficient accuracy. However, all of the accelerated ADMM methods outperform ADMM in convergence of primal residual. We can see in Fig.~\ref{dual} that the convergence rate of dual residual is better compared to convergence rate of primal residual for ADMM, notably achieving accuracy up to $10^{-3}$ after 200 iterations. However, like in convergence of primal residual, all accelerated ADMM methods outperform ADMM. 

We note that restarting the algorithm in every 50 iterations makes the convergence of the primal residual and dual residual significantly faster compared to ADMM and other accelerated ADMM with different restart schemes. Also, we can see that restarting every 10 or 20 iterations makes the convergence rate of primal residual slower compared to accelerated ADMM without restarting the algorithm in Fig.~\ref{primal}. These two restart schemes does not have a significant effect of convergence rate of dual residual compared to accelerated ADMM without restart. Even though some of the restart methods does not perform as well as restarting in every 50 iterations, all of the accelerated variants outperforms ADMM. After 200 iterations with the accelerated ADMM, the minimum (upper-bound on)  induced $L_2$-norm from the  input $d$ to the output output $y$  is $0.29$. Also, we note that the convergence of the accelerated ADMM becomes irregular when the residuals are below $10^{-6}$, and this phenomenon can be explained by the solver precision of solving the local SDP problems, which are usually set to $10^{-6}$ for practical reasons.

\section{Conclusions and Future Work}\label{sec:conclude}

We proposed a method for compositional analysis  of large scale  hybrid systems defined over finite alphabets, for which an underlying interconnection topology is given. For such systems, we decompose the global analysis problem into a number of smaller local analysis problems using dissipativity theory. Furthermore, we proposed a distributed optimization method with smoothing techniques, which enables to employ accelerated ADMM. Numerical results show that the accelerated ADMM method with different restart methods significantly improves the convergence rate compared to standard ADMM. 

In the examples studied in this paper, we considered linear continuous dynamics. The generalization to polynomial continuous dynamics can be formulated based on sum-of-squares optimization~(\cite{Par00}). Moreover, one interesting analysis problem for future research is compositional safety verification. In this respect, in~\cite{Sloth2012}, a method is brought forward based on compositional barrier certificates for continuous systems and the dual decomposition method was used for implementation. Finally, distributed synthesis of control laws for large-scale hybrid systems can be studied using the dissipativity framework presented in this study. 
  
\appendix

\section{Proof of Theorem 1} \label{app1}
 Inequality~\eqref{inq1stab} implies that $V$ is positive definite and $V(0,q)=0,~\forall q \in \mathscr{Q}$. We define the time intervals $T_k = (t_k,t_{k+1})$ where the continuous dynamics follows $\mathcal{C}$ in \eqref{eq:HS} with $p \in \mathscr{P}$. Similarly, $t^{+}_k$  corresponds to the discrete jump instant  at $t_k$. By continuity, we can re-write~\eqref{inq2stab} as
\begin{multline} \label{dddfd}
   \left( \frac{\partial V(x,q)}{\partial x} \right)^T f(x,0;p) \le -\rho_p,\quad \forall x \in \mathbb{R}^n, \\ \forall q \in \mathscr{Q},~~\forall p \in \mathscr{P},
 \end{multline}
 where $\rho_p,~p \in \mathscr{P}$, is a small positive number. Moreover, differentiating  $V$ with respect to time and integrating it from $t_0$ to $t_n$ and noting that the discrete jumps happen on sets of measure zero, yields
\begin{multline} \label{ssddf}
\int_{t_0}^{t_n} \frac{dV}{dt}~dt =    \sum_{k=0}^{n-1} \int_{t_k}^{t_{k+1}}   \left( \frac{\partial V(x,q)}{\partial x} \right)^T f(x,0;p_k)~dt \\+ \sum_{k=0}^{n-1} \left(V\big(x(t_{k+1}),q^{+}(t_{k+1})\big) - V\big(x(t_{k+1}),q(t_{k+1})\big) \right).
\end{multline}
From~\eqref{popos}, we infer
$$
\sum_{k=0}^{n-1} \left(V\big(x(t_{k+1}),q^{+}(t_{k+1})\big) - V\big(x(t_{k+1}),q(t_{k+1})\big) \right) \le 0,
$$
since it is the finite sum of non-positive terms. Hence, 
$$
\int_{t_0}^{t_n} \frac{dV}{dt}~dt \le \sum_{k=0}^{n-1} \int_{t_k}^{t_{k+1}}   \left( \frac{\partial V(x,q)}{\partial x} \right)^T f(x,0;p_k)~dt
$$
 Using~\eqref{dddfd}, the right-hand side of above inequality  satisfies
\begin{align} \label{scc}
\sum_{k=0}^{n-1} \int_{t_k}^{t_{k+1}}   \left( \frac{\partial V(x,q)}{\partial x} \right)^T f(x,0;p_k)~dt \nonumber \\
\le -\int_{t_0}^{t_n} \rho~dt = -\rho (t_n - t_0),
\end{align}
where $\rho = \min_{p \in \mathscr{P}} \rho_p$. That is,
$$
\int_{t_0}^{t_n} \frac{dV}{dt}~dt = V(x(t_n),q(t_n)) - V(x_0,q_0) \le -\rho (t_n - t_0)
$$
Re-organizing the terms gives 
$$
V(x(t_n),q(t_n)) \le  V(x_0,q_0) +\rho t_0 -\rho t_n.
$$
From~\eqref{inq1stab}, we know that $V(x(t_n),q(t_n))>0$; therefore, there exists a $t_n \ge \frac{\rho_0t_0+V(x_0,q_0)}{\rho}$ such that $V\left(x(t_n),q(t_n)\right)=0$. Additionally, since~\eqref{popos} holds,  we have $V(x,q^{+}(t)) \le V(x,q(t))$. Thus, $V(x,q(t))=0$ for all $q \in \mathscr{Q}$ which implies  $x(t)=0$ for all $t \ge t_n$.

\section{Proof of Proposition 3} \label{app3}
Multiplying the inequality in~\eqref{sddccxcxglobal} from left by $\begin{bmatrix} y \\ d \end{bmatrix}^T$ and right by $\begin{bmatrix} y \\ d \end{bmatrix}$, we obtain
\begin{equation} \label{op11}
\sum_{i=1}^N \begin{bmatrix} w_i \\ y_i \end{bmatrix}^T S_i \begin{bmatrix} w_i \\ y_i \end{bmatrix} - \begin{bmatrix} d \\ z \end{bmatrix}^T S \begin{bmatrix} d \\ z \end{bmatrix} \le 0.
\end{equation}
Similarly, multiplying the inequality in~\eqref{sddccxcxglobald} from left by $\begin{bmatrix} p \\ \mu \end{bmatrix}^T$ and right by $\begin{bmatrix} p \\ \mu \end{bmatrix}$ gives
\begin{equation} \label{op22}
\sum_{i=1}^N \begin{bmatrix} u_i \\ p_i \end{bmatrix}^T R_i \begin{bmatrix} u_i \\ p_i \end{bmatrix} - \begin{bmatrix} \mu \\ \zeta \end{bmatrix}^T R \begin{bmatrix} \mu \\ \zeta \end{bmatrix} \le 0.
\end{equation}
Moreover, because $(S_i,R_i) \in \mathcal{L}_i$, the exists storage functions $V_i(x,q),~i=1,2,\ldots,N$ such that
$$
 \left( \frac{\partial V_i(x_i,q_i)}{\partial x_i} \right)^T f_i(x_i,w_i;p_i) - \begin{bmatrix} w_i \\ y_i \end{bmatrix}^T S_i \begin{bmatrix} w_i \\ y_i \end{bmatrix} \le 0,
$$
and
$$
  {V}(x_i,q^{+}_i)- V(x_i,q_i)  - \begin{bmatrix} u_i \\ p_i \end{bmatrix}^T R_i \begin{bmatrix} u_i \\ p_i \end{bmatrix} \le 0.
$$
If we sum over $i=1$ to $N$ the above dissipation inequalities and use \eqref{op11} and \eqref{op22}, we infer
$$
\sum_{i=1}^N  \left( \frac{\partial V_i(x_i,q_i)}{\partial x_i} \right)^T f_i(x_i,w_i;p_i) \le \begin{bmatrix} d \\ z \end{bmatrix}^T S \begin{bmatrix} d \\ z \end{bmatrix},
$$
and 
$$
\sum_{i=1}^N \left(   {V}(x_i,q^{+}_i)- V(x_i,q_i)  \right) \le \begin{bmatrix} \mu \\ \zeta \end{bmatrix}^T R \begin{bmatrix} \mu \\ \zeta \end{bmatrix},
$$
which implies that the overall system is dissipative with the storage function $V(x,q) = \sum_{i=1}^N V_i(x_i,q_i)$.

\bibliographystyle{natbib}  

\bibliography{ifacconf}

\end{document}

%% file: ADHS_hybrid_primal.tex
%
%
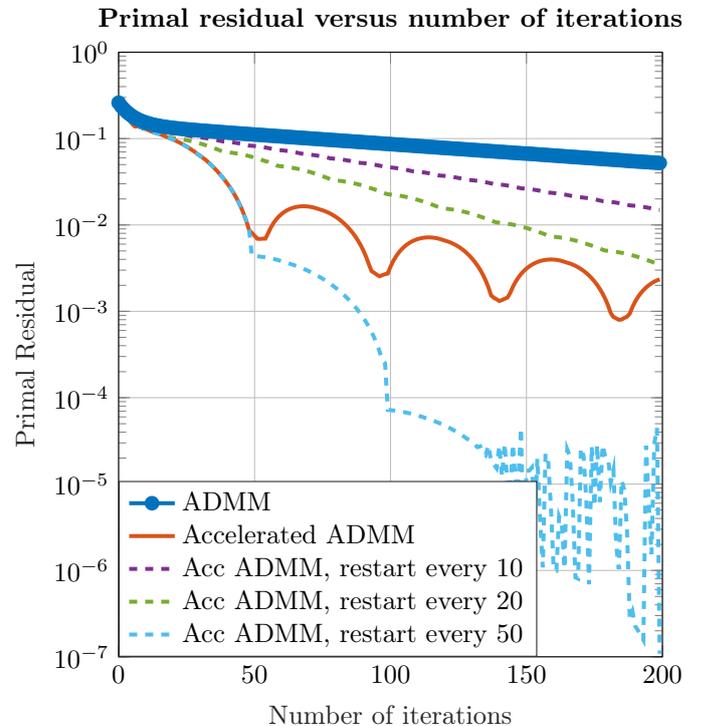
\begin{figure}[H]
\definecolor{mycolor1}{rgb}{0.00000,0.44700,0.74100}%
\definecolor{mycolor2}{rgb}{0.85000,0.32500,0.09800}%
\definecolor{mycolor3}{rgb}{0.92900,0.69400,0.12500}%
\definecolor{mycolor4}{rgb}{0.49400,0.18400,0.55600}%
\definecolor{mycolor5}{rgb}{0.46600,0.67400,0.18800}%
\definecolor{mycolor6}{rgb}{0.30100,0.74500,0.93300}%
\begin{tikzpicture}

\begin{axis}[%
width=2.817in,
height=3.156in,
at={(2.167in,0.898in)},
scale only axis,
xmajorgrids=true,
ymajorgrids=true,
xmin=0,
xmax=200,
xlabel style={font=\color{white!15!black}},
xlabel={Number of iterations},
ymode=log,
ymin=1e-07,
ymax=1,
yminorticks=true,
legend style={at={(0,0)},anchor=south west,font=\fontsize{10}{10}\selectfont},
ylabel style={font=\color{white!15!black}},
ylabel={Primal Residual},
axis background/.style={fill=white},
title style={font=\bfseries},
title={Primal residual versus number of iterations},
legend style={legend cell align=left, align=left, draw=white!15!black}
]
\addplot [color=mycolor1, line width=1.5pt, draw=none, mark=*, mark options={solid, mycolor1}]
  table[row sep=crcr]{%
0	0.260313969825423\\
1	0.236015452636075\\
2	0.218289662593266\\
3	0.203524355485969\\
4	0.191783316797702\\
5	0.182172643212716\\
6	0.174202939390984\\
7	0.167556315725182\\
8	0.161979671830751\\
9	0.157276130478697\\
10	0.153286274919383\\
11	0.149878311808575\\
12	0.146943030103828\\
13	0.144400827250925\\
14	0.14217806873792\\
15	0.140218490028999\\
16	0.138475539282087\\
17	0.136909483414083\\
18	0.135493291718239\\
19	0.134199756355039\\
20	0.133008373005484\\
21	0.131903749461001\\
22	0.130871237294326\\
23	0.129896302140327\\
24	0.128973377970702\\
25	0.128095618013391\\
26	0.127252539760021\\
27	0.126439943673101\\
28	0.125654612806773\\
29	0.124892623957058\\
30	0.124150901334686\\
31	0.123426886737012\\
32	0.122718463199903\\
33	0.122023855121755\\
34	0.121341551180376\\
35	0.120670283197698\\
36	0.120008967491605\\
37	0.119356599028146\\
38	0.118713402432276\\
39	0.118077213279557\\
40	0.117448128085579\\
41	0.116825150375639\\
42	0.116208089142106\\
43	0.115596672510509\\
44	0.114990618217144\\
45	0.114389416480709\\
46	0.11379290131636\\
47	0.113200864831602\\
48	0.11261312475931\\
49	0.112029517717568\\
50	0.111449901737511\\
51	0.110874151082874\\
52	0.110302151552451\\
53	0.109733925852329\\
54	0.109169358141937\\
55	0.108608373650722\\
56	0.108050952677395\\
57	0.10749676707612\\
58	0.106945769986221\\
59	0.106397917832834\\
60	0.105853278297278\\
61	0.105311669565753\\
62	0.104773039731384\\
63	0.104237367268268\\
64	0.103704624491567\\
65	0.103174789847218\\
66	0.102647838492517\\
67	0.102123744613304\\
68	0.10160248404679\\
69	0.101084015551853\\
70	0.100568276194281\\
71	0.100055304671836\\
72	0.0995450821703084\\
73	0.099037586793132\\
74	0.0985327871202262\\
75	0.0980306747812473\\
76	0.0975312285646544\\
77	0.0970344259402542\\
78	0.0965402516505325\\
79	0.0960486857731693\\
80	0.0955597097565434\\
81	0.0950733089893157\\
82	0.0945894679294592\\
83	0.0941081683719861\\
84	0.0936293940303189\\
85	0.0931531291442077\\
86	0.0926793525577737\\
87	0.0922080461356422\\
88	0.0917392171412481\\
89	0.0912728374759408\\
90	0.0908088888824407\\
91	0.0903475209060781\\
92	0.08988795988862\\
93	0.089430660717989\\
94	0.0889757987889271\\
95	0.0885233449854249\\
96	0.0880732717015428\\
97	0.0876255562572945\\
98	0.0871801777229324\\
99	0.0867371153572788\\
100	0.0862962880815061\\
101	0.0858577539036703\\
102	0.0854214981013733\\
103	0.0849874613281066\\
104	0.0845556774787646\\
105	0.0841261368000296\\
106	0.0836988240860022\\
107	0.0832735555556824\\
108	0.0828504070866917\\
109	0.0824294845636598\\
110	0.0820107681025467\\
111	0.0815942619524309\\
112	0.0811799350940664\\
113	0.0807677731768252\\
114	0.0803577579023713\\
115	0.0799498706237464\\
116	0.0795440975254184\\
117	0.0791404243177183\\
118	0.0787388347014919\\
119	0.0783393151021204\\
120	0.0779418570168509\\
121	0.0775464480098344\\
122	0.0771530763434052\\
123	0.0767617293506956\\
124	0.0763723957790146\\
125	0.0759850647663614\\
126	0.0755997250716803\\
127	0.0752163662423399\\
128	0.0748349783730647\\
129	0.074455561130787\\
130	0.0740780609848695\\
131	0.0737024822644163\\
132	0.0733288445711346\\
133	0.0729571280472878\\
134	0.0725873122573443\\
135	0.0722193861158459\\
136	0.0718533467365629\\
137	0.0714891841785413\\
138	0.071126889251151\\
139	0.0707664503264331\\
140	0.0704078573378114\\
141	0.0700511003340668\\
142	0.0696961687080374\\
143	0.0693430519438491\\
144	0.0689917410529348\\
145	0.0686422286457243\\
146	0.0682945015949081\\
147	0.067948546898822\\
148	0.0676043562181435\\
149	0.0672619199284588\\
150	0.0671605888027058\\
151	0.0667785740251475\\
152	0.0664072599426714\\
153	0.0660446547755804\\
154	0.0656892613085534\\
155	0.0653399099798931\\
156	0.0649956801280955\\
157	0.0646558448175823\\
158	0.0643198305350561\\
159	0.0639871807806653\\
160	0.0636575272876298\\
161	0.063330621251973\\
162	0.0630061992297559\\
163	0.0626840926366418\\
164	0.0623641480105512\\
165	0.0620462315922828\\
166	0.0617302894089254\\
167	0.0614162312714907\\
168	0.0611039886161143\\
169	0.0607935073977536\\
170	0.060484743811193\\
171	0.0601776620426625\\
172	0.0598722313272493\\
173	0.059568425051466\\
174	0.0592662204871964\\
175	0.0589655940143739\\
176	0.0586665333412765\\
177	0.0583690235033495\\
178	0.0580730592937039\\
179	0.0577786078235308\\
180	0.057485669986717\\
181	0.0571942353447004\\
182	0.0569042934872134\\
183	0.0566158344430377\\
184	0.0563288476241974\\
185	0.0560433258824928\\
186	0.0557592580610918\\
187	0.055275200160103\\
188	0.0550317528737935\\
189	0.0547818182568542\\
190	0.0545271893411131\\
191	0.0542691050655619\\
192	0.0540085991730938\\
193	0.0537464675933363\\
194	0.0534833336998777\\
195	0.0532196887254963\\
196	0.0529559387665843\\
197	0.052692366643081\\
198	0.0524292236914604\\
199	0.0521666573660947\\
};
\addlegendentry{ADMM}

\addplot [color=mycolor2, line width=1.5pt]
  table[row sep=crcr]{%
0	0.260313969825423\\
1	0.229867301044826\\
2	0.203716220055454\\
3	0.18038626129511\\
4	0.161778184835114\\
5	0.14750692787763\\
6	0.137963044527752\\
7	0.13833801739002\\
8	0.136559218975332\\
9	0.133329066212017\\
10	0.129316771459057\\
11	0.125064406489225\\
12	0.120966824220796\\
13	0.117646552719208\\
14	0.115549653086002\\
15	0.112986417629386\\
16	0.110010111081072\\
17	0.106725825220738\\
18	0.103238296583957\\
19	0.100404758144793\\
20	0.0975257978797006\\
21	0.0944579045984338\\
22	0.0912158856986584\\
23	0.0878412367301605\\
24	0.0843772728709166\\
25	0.0808634806556223\\
26	0.0773254383603127\\
27	0.0737850912062463\\
28	0.0702560427403155\\
29	0.0662711587506805\\
30	0.0631479271825446\\
31	0.0597872446924505\\
32	0.0568166641598311\\
33	0.0531235971476436\\
34	0.0496603618687457\\
35	0.0462253737802097\\
36	0.0428445673952794\\
37	0.0395376661695471\\
38	0.0363131278999792\\
39	0.0331697571607362\\
40	0.0301042755469394\\
41	0.0271436253636023\\
42	0.0242644250684383\\
43	0.0214610304232167\\
44	0.0187297617591223\\
45	0.0160760626145048\\
46	0.0135060526805539\\
47	0.0110269313990558\\
48	0.00870348158633383\\
49	0.00810536705427206\\
50	0.00750932997900761\\
51	0.00691726411261251\\
52	0.00689146196039248\\
53	0.00693633544465966\\
54	0.00694658872954512\\
55	0.00803764619685415\\
56	0.00935605532444495\\
57	0.0105563351268261\\
58	0.0116388209198395\\
59	0.0126056794391894\\
60	0.0134579689201216\\
61	0.0141954234099835\\
62	0.0148223728043581\\
63	0.0153423017022877\\
64	0.0157588859796415\\
65	0.0160758290354571\\
66	0.0162969066067822\\
67	0.0164258019719681\\
68	0.0164662228465183\\
69	0.0164220496237684\\
70	0.0162972252990459\\
71	0.0160957083437203\\
72	0.0158216983171554\\
73	0.0154795965891361\\
74	0.015116006516233\\
75	0.0147890232104582\\
76	0.0144014071190636\\
77	0.0139574693707734\\
78	0.0134615561947338\\
79	0.0129181009533325\\
80	0.012331569996514\\
81	0.0117064597520183\\
82	0.0110461920537334\\
83	0.0103560224211391\\
84	0.00964211473281285\\
85	0.008907601699226\\
86	0.00815677708590862\\
87	0.00739345734106084\\
88	0.00662163107076452\\
89	0.00584498065512346\\
90	0.00506672849746892\\
91	0.00429249378231439\\
92	0.00352520456956384\\
93	0.00292583071055319\\
94	0.00279201594441159\\
95	0.0026500556448702\\
96	0.0025432014482179\\
97	0.0026248296948963\\
98	0.00269076041356325\\
99	0.00274115401182723\\
100	0.00328674525057776\\
101	0.00382111837188853\\
102	0.0043185124596582\\
103	0.00477767913224642\\
104	0.00519766539320155\\
105	0.00557772913970733\\
106	0.00591734029510072\\
107	0.00621617749812453\\
108	0.00647411909654966\\
109	0.00669123014844281\\
110	0.00686775609145485\\
111	0.00700411001898693\\
112	0.00710086902958118\\
113	0.00715876658322817\\
114	0.00717868672336885\\
115	0.00716166225483411\\
116	0.00710886851155945\\
117	0.00702161324089723\\
118	0.00690133434879037\\
119	0.00674958638241463\\
120	0.00659298742775222\\
121	0.00644875151616081\\
122	0.00627606447398782\\
123	0.00607658233384362\\
124	0.0058521791293033\\
125	0.0056041917891832\\
126	0.00533500158288214\\
127	0.00504638946398711\\
128	0.00474027255751225\\
129	0.00441862317535086\\
130	0.00408344743960321\\
131	0.00373676032149645\\
132	0.00338057880300235\\
133	0.00301689960469626\\
134	0.00264771501528333\\
135	0.00227499089132643\\
136	0.00190063534436136\\
137	0.00152654221734538\\
138	0.00143757199459105\\
139	0.00138001397028736\\
140	0.00131714787457504\\
141	0.00135140207896156\\
142	0.00139965419743504\\
143	0.00143999290208141\\
144	0.00162897640852744\\
145	0.0019242552191811\\
146	0.00220246246975847\\
147	0.00246264347772351\\
148	0.00270392411097375\\
149	0.0029256411977085\\
150	0.00312716072554342\\
151	0.00330797089566104\\
152	0.00346767528927855\\
153	0.00360599224224362\\
154	0.00372274996020633\\
155	0.00381788709978925\\
156	0.00389144781671236\\
157	0.00394357922237987\\
158	0.00397452400058507\\
159	0.00398460034219234\\
160	0.00397431987686642\\
161	0.00394416350907328\\
162	0.00389472902675603\\
163	0.00382664605042576\\
164	0.00374074018008302\\
165	0.00365679441126349\\
166	0.00357519509679199\\
167	0.00347734150955309\\
168	0.0033640842975271\\
169	0.00323632432911419\\
170	0.00309503027667163\\
171	0.00294120606210757\\
172	0.00277589066388961\\
173	0.00260015944568642\\
174	0.00241512347915052\\
175	0.0022218973282182\\
176	0.00202162874237588\\
177	0.00181545729142252\\
178	0.0016045388166574\\
179	0.00139003042499036\\
180	0.00117306166621395\\
181	0.000954771348824379\\
182	0.000859123756774771\\
183	0.000829816237773562\\
184	0.000796989858454111\\
185	0.000805315052020916\\
186	0.00083906956472718\\
187	0.000867962764358798\\
188	0.000942268574318284\\
189	0.00112568337979335\\
190	0.00129966297434253\\
191	0.00146353208443251\\
192	0.00161671040871586\\
193	0.00175858259981787\\
194	0.00188872309991717\\
195	0.00200666201185235\\
196	0.00211208175581529\\
197	0.00220472279485645\\
198	0.00228441780631888\\
199	0.00235096046647076\\
};
\addlegendentry{Accelerated ADMM}

\addplot [color=mycolor4, dashed, line width=1.5pt]
  table[row sep=crcr]{%
0	0.260313969825423\\
1	0.229867301044826\\
2	0.203716220055454\\
3	0.18038626129511\\
4	0.161778184835114\\
5	0.14750692787763\\
6	0.137963044527752\\
7	0.13833801739002\\
8	0.136559218975332\\
9	0.1305366795155\\
10	0.128877524753928\\
11	0.127101948201447\\
12	0.126237501817086\\
13	0.125610719512308\\
14	0.124612713601209\\
15	0.123242151430346\\
16	0.121559507746979\\
17	0.119609069546593\\
18	0.117443101888908\\
19	0.116709579855055\\
20	0.115883824709361\\
21	0.114903943246886\\
22	0.113779562563668\\
23	0.112518354006579\\
24	0.111126026134675\\
25	0.109606974952803\\
26	0.107965582990118\\
27	0.106206047606263\\
28	0.104329099050227\\
29	0.103771755924831\\
30	0.103096580984379\\
31	0.10228091554694\\
32	0.101329736736934\\
33	0.100247515757059\\
34	0.0990386365625559\\
35	0.0977078036339112\\
36	0.0962595542635334\\
37	0.0946986580943291\\
38	0.0930300047428902\\
39	0.0925445043728549\\
40	0.0919420127403935\\
41	0.0912148337903655\\
42	0.09036741967009\\
43	0.0894034579243521\\
44	0.0883255989352006\\
45	0.0871400640645717\\
46	0.0858502067653553\\
47	0.0844597340031886\\
48	0.082973359425113\\
49	0.0825450435790725\\
50	0.0820082310188738\\
51	0.0813603159510114\\
52	0.0806052158132749\\
53	0.0797465163332236\\
54	0.0787876880930188\\
55	0.077732218531899\\
56	0.0765836762003371\\
57	0.0753457515922995\\
58	0.0740222420571232\\
59	0.0736436416804054\\
60	0.0731652752147275\\
61	0.0725879304038969\\
62	0.0719150674158775\\
63	0.0711498965493017\\
64	0.0702955105387779\\
65	0.0693550063229646\\
66	0.0683315431479099\\
67	0.0672283822905479\\
68	0.0664718337907349\\
69	0.0658814667282182\\
70	0.0654182786604304\\
71	0.0648695451804772\\
72	0.0642407318787666\\
73	0.0635357734730281\\
74	0.0627575065547779\\
75	0.0619080261939425\\
76	0.0609890151068293\\
77	0.0600020723087738\\
78	0.0589488945718176\\
79	0.0586504762537843\\
80	0.0582712439420099\\
81	0.0578130106201125\\
82	0.0572783852360757\\
83	0.0566698218739174\\
84	0.0559897412633643\\
85	0.0552406099939319\\
86	0.054074819975961\\
87	0.0534539414406299\\
88	0.0525762818718951\\
89	0.052311450068202\\
90	0.0519739581522921\\
91	0.0515661466535189\\
92	0.0510903118643225\\
93	0.0505485781150654\\
94	0.0499430518179677\\
95	0.0492758886317085\\
96	0.0485493332211655\\
97	0.0477657062381125\\
98	0.0469274576569951\\
99	0.0466894670830827\\
100	0.0463863166439187\\
101	0.0460204320004106\\
102	0.0455939961090734\\
103	0.0451090028343626\\
104	0.044567371326527\\
105	0.0439710998431041\\
106	0.0433221156958229\\
107	0.0426224935303682\\
108	0.0418743534221151\\
109	0.0416622832687329\\
110	0.0413919060180366\\
111	0.0410655538407934\\
112	0.0406851669698075\\
113	0.0402525167883276\\
114	0.0397693282652716\\
115	0.0392373214476814\\
116	0.0386583207407267\\
117	0.0380341524815125\\
118	0.0373666861510885\\
119	0.0371775842963953\\
120	0.0369363363441197\\
121	0.0366451418136248\\
122	0.0363057236179009\\
123	0.0359196858682021\\
124	0.0354885475512959\\
125	0.0350139356012741\\
126	0.0344984581531144\\
127	0.0339418538076794\\
128	0.0333463727955409\\
129	0.0331775895274304\\
130	0.0329621934780248\\
131	0.0327022341466225\\
132	0.0323992543290924\\
133	0.0320546716421108\\
134	0.0316698573350059\\
135	0.0312461740012578\\
136	0.0307869866598473\\
137	0.0302911790169866\\
138	0.029759720213916\\
139	0.0296088060089684\\
140	0.0294164440970736\\
141	0.0291842837740938\\
142	0.02891377979419\\
143	0.0286061308630309\\
144	0.028262583884935\\
145	0.0278843735955111\\
146	0.0274727575984234\\
147	0.0270290373362175\\
148	0.0265545569592148\\
149	0.0264203002017802\\
150	0.0262488594909324\\
151	0.0260419411030431\\
152	0.0258007554918469\\
153	0.0255264110805961\\
154	0.0252199954279201\\
155	0.0248826231156179\\
156	0.0245154183795949\\
157	0.0241195435883709\\
158	0.0236962002239115\\
159	0.0235762995822091\\
160	0.0234232535820085\\
161	0.0232385380243097\\
162	0.0230232493779676\\
163	0.0227783902868426\\
164	0.0225049379139781\\
165	0.0222038739992821\\
166	0.0218762025341862\\
167	0.0215229582788429\\
168	0.0211452119800773\\
169	0.0210383729650143\\
170	0.0209018438321023\\
171	0.0207370519870258\\
172	0.0205449747788907\\
173	0.0203265067781703\\
174	0.0200825188573033\\
175	0.0198138875484672\\
176	0.0195215098184049\\
177	0.0192063091453443\\
178	0.018869242104921\\
179	0.0187739223625669\\
180	0.0186520825955336\\
181	0.0185050234043586\\
182	0.0183336182820016\\
183	0.0181387043750274\\
184	0.0179210262097483\\
185	0.0176813299878686\\
186	0.0174204293801939\\
187	0.0171391504999091\\
188	0.0168383510749296\\
189	0.0167532966284438\\
190	0.0166445604924566\\
191	0.0165133210057506\\
192	0.0163603580104838\\
193	0.01618638464207\\
194	0.0159920940287488\\
195	0.0157781836946315\\
196	0.0155453667086411\\
197	0.0152943780958993\\
198	0.0150259780019831\\
199	0.0149500970375909\\
};
\addlegendentry{Acc ADMM, restart every 10}

\addplot [color=mycolor5, dashed, line width=1.5pt]
  table[row sep=crcr]{%
0	0.260313969825423\\
1	0.229867301044826\\
2	0.203716220055454\\
3	0.18038626129511\\
4	0.161778184835114\\
5	0.14750692787763\\
6	0.137963044527752\\
7	0.13833801739002\\
8	0.136559218975332\\
9	0.133329066212017\\
10	0.129316771459057\\
11	0.125064406489225\\
12	0.120966824220796\\
13	0.117646552719208\\
14	0.115549653086002\\
15	0.112986417629386\\
16	0.110010111081072\\
17	0.106725825220738\\
18	0.103238296583957\\
19	0.102546592366774\\
20	0.101755121849056\\
21	0.10082790872302\\
22	0.0997777206075678\\
23	0.0986147236037034\\
24	0.0973461617133195\\
25	0.0959772047609083\\
26	0.0945114313467073\\
27	0.0929514026060256\\
28	0.0912991355978977\\
29	0.0895564984161562\\
30	0.0877232848585018\\
31	0.0858025652160656\\
32	0.0837994330599748\\
33	0.0817168780791873\\
34	0.0797488666568303\\
35	0.0777606016102458\\
36	0.0757012605932569\\
37	0.0735761077248273\\
38	0.0713907272116431\\
39	0.0699261970887839\\
40	0.0694698171029882\\
41	0.0689195518526988\\
42	0.0682788701511034\\
43	0.0675508789094361\\
44	0.0671285624826129\\
45	0.0659672108480248\\
46	0.0649290726735279\\
47	0.0638363345390955\\
48	0.0626914711536485\\
49	0.0614938474371348\\
50	0.0602414565685172\\
51	0.0589324656240893\\
52	0.0575660681745408\\
53	0.0561429237284983\\
54	0.0547633845861737\\
55	0.0537811867429564\\
56	0.05204684176292\\
57	0.0505170300334579\\
58	0.0489723537755156\\
59	0.0480760160788534\\
60	0.047751425230971\\
61	0.0473635391717232\\
62	0.046915410915133\\
63	0.0464094461528262\\
64	0.045847601283118\\
65	0.0452315847108842\\
66	0.0445629741102348\\
67	0.0438433877266565\\
68	0.0430744901369587\\
69	0.0422581523395003\\
70	0.0413964292210647\\
71	0.0404915636362962\\
72	0.0395459720366829\\
73	0.0385621987071436\\
74	0.0376293067421034\\
75	0.0366953429747813\\
76	0.0357276631970608\\
77	0.0347286024995562\\
78	0.0336982704791416\\
79	0.0330014337547244\\
80	0.0327871867095603\\
81	0.0325286118926\\
82	0.0322272468603835\\
83	0.0318845021784018\\
84	0.0315017401224333\\
85	0.0310805982959401\\
86	0.0306243608645957\\
87	0.0301304942540263\\
88	0.0296015821579922\\
89	0.0290396057278255\\
90	0.0284463805939735\\
91	0.0278236467720832\\
92	0.0271731812700922\\
93	0.0264967820408478\\
94	0.0258603997243112\\
95	0.0252185741305114\\
96	0.024553536853172\\
97	0.0238671267868089\\
98	0.0231613927805515\\
99	0.0226739451612947\\
100	0.0225271801277006\\
101	0.0223499219731079\\
102	0.0221431997051456\\
103	0.0219079655151249\\
104	0.0216451557621372\\
105	0.021355721473276\\
106	0.021040640026983\\
107	0.0207009215485774\\
108	0.0203376102864041\\
109	0.0199517823356241\\
110	0.019544549110332\\
111	0.0191170523813529\\
112	0.0186704638712805\\
113	0.0182059842604435\\
114	0.0177722210443056\\
115	0.0173309456612463\\
116	0.0168737649627576\\
117	0.0164018903011197\\
118	0.0159165196336672\\
119	0.015580229895161\\
120	0.0154791953392271\\
121	0.0153572264520507\\
122	0.0152150428447364\\
123	0.0150533045244476\\
124	0.0148726561696403\\
125	0.0146737472636867\\
126	0.0144572436504951\\
127	0.0142238316685533\\
128	0.0139742215625152\\
129	0.0137091467976144\\
130	0.013429366022346\\
131	0.0131356601314311\\
132	0.0128288326239155\\
133	0.012509707848738\\
134	0.0122130241590606\\
135	0.0119097792175393\\
136	0.0115955870528706\\
137	0.0112712685086186\\
138	0.0109376569390485\\
139	0.0107055768184574\\
140	0.0106361322627378\\
141	0.0105523052290722\\
142	0.0104545917184142\\
143	0.0103434462582473\\
144	0.0102193118028999\\
145	0.0100826345187102\\
146	0.00993387123157985\\
147	0.00977349278850265\\
148	0.00960198620238922\\
149	0.00941985498193226\\
150	0.00922761909594438\\
151	0.00902581482676268\\
152	0.00881499401245467\\
153	0.00859572272258047\\
154	0.00839242779282207\\
155	0.0081840212072759\\
156	0.00796809297121676\\
157	0.00774520719808874\\
158	0.00751593362708292\\
159	0.00735609766913742\\
160	0.00730836545093694\\
161	0.00725075199604724\\
162	0.00718359898194271\\
163	0.00710721909490636\\
164	0.00702191687298404\\
165	0.00692799879881118\\
166	0.00682577866116307\\
167	0.00671557892621044\\
168	0.00659773008536035\\
169	0.0064725919946641\\
170	0.00634051967129173\\
171	0.00620187905176672\\
172	0.0060570455106862\\
173	0.00590640248151481\\
174	0.00576693845897591\\
175	0.00562371143528684\\
176	0.00547545583722009\\
177	0.00532223256298836\\
178	0.00516447328646191\\
179	0.00505465604964267\\
180	0.00502185361884155\\
181	0.00498223670031073\\
182	0.0049360675911887\\
183	0.00488356629863113\\
184	0.00482494281523394\\
185	0.00476040579480269\\
186	0.00469016902110425\\
187	0.00461445175515829\\
188	0.00453348379331536\\
189	0.00444752853093927\\
190	0.00435676669371526\\
191	0.00426147389207948\\
192	0.0041619223325064\\
193	0.00405838276147286\\
194	0.00396266662964662\\
195	0.00386426055529876\\
196	0.00376229280986939\\
197	0.00365704055793725\\
198	0.00354876718541836\\
199	0.00347311069551542\\
};
\addlegendentry{Acc ADMM, restart every 20}

\addplot [color=mycolor6, dashed, line width=1.5pt]
  table[row sep=crcr]{%
0	0.260313969825423\\
1	0.229867301044826\\
2	0.203716220055454\\
3	0.18038626129511\\
4	0.161778184835114\\
5	0.14750692787763\\
6	0.137963044527752\\
7	0.13833801739002\\
8	0.136559218975332\\
9	0.133329066212017\\
10	0.129316771459057\\
11	0.125064406489225\\
12	0.120966824220796\\
13	0.117646552719208\\
14	0.115549653086002\\
15	0.112986417629386\\
16	0.110010111081072\\
17	0.106725825220738\\
18	0.103238296583957\\
19	0.100404758144793\\
20	0.0975257978797006\\
21	0.0944579045984338\\
22	0.0912158856986584\\
23	0.0878412367301605\\
24	0.0843772728709166\\
25	0.0808634806556223\\
26	0.0773254383603127\\
27	0.0737850912062463\\
28	0.0702560427403155\\
29	0.0662711587506805\\
30	0.0631479271825446\\
31	0.0597872446924505\\
32	0.0568166641598311\\
33	0.0531235971476436\\
34	0.0496603618687457\\
35	0.0462253737802097\\
36	0.0428445673952794\\
37	0.0395376661695471\\
38	0.0363131278999792\\
39	0.0331697571607362\\
40	0.0301042755469394\\
41	0.0271436253636023\\
42	0.0242644250684383\\
43	0.0214610304232167\\
44	0.0187297617591223\\
45	0.0160760626145048\\
46	0.0135060526805539\\
47	0.0110269313990558\\
48	0.00870348158633383\\
49	0.00438534085282366\\
50	0.00436477233394944\\
51	0.00433797556266274\\
52	0.00430463562388989\\
53	0.00426461444664167\\
54	0.0042179513419528\\
55	0.0041648365276743\\
56	0.00410557786467253\\
57	0.00404055833633332\\
58	0.0039701993451482\\
59	0.0038949291292561\\
60	0.00381516009014588\\
61	0.00373127542196205\\
62	0.00364362449137678\\
63	0.00355252278408413\\
64	0.00346870900890206\\
65	0.00338270495415927\\
66	0.00329320144841117\\
67	0.0032009667608674\\
68	0.00310612295130508\\
69	0.00300889169322455\\
70	0.00290952660171273\\
71	0.00280825927920775\\
72	0.00270535249413167\\
73	0.00260104270499262\\
74	0.00249557548125384\\
75	0.00238919681728245\\
76	0.00228215156184625\\
77	0.00217467689041551\\
78	0.002067014748382\\
79	0.00195940601623645\\
80	0.00185207908221763\\
81	0.00174526672004044\\
82	0.00163932347927798\\
83	0.00153401853269247\\
84	0.00143005666917906\\
85	0.00132751173586726\\
86	0.00122654477131762\\
87	0.00112733507230017\\
88	0.00103006885455958\\
89	0.000934917155860582\\
90	0.000842054818605831\\
91	0.000751634595546835\\
92	0.000663804863239359\\
93	0.000578703604873387\\
94	0.000496465996009554\\
95	0.000417194331295501\\
96	0.00034099125209397\\
97	0.000268750126038997\\
98	0.000250308246840691\\
99	7.22522212660559e-05\\
100	7.17820154572113e-05\\
101	7.12140783478144e-05\\
102	7.05526794321787e-05\\
103	6.98002810137766e-05\\
104	6.89607309657947e-05\\
105	6.80361963887322e-05\\
106	6.70304737058203e-05\\
107	6.59466635529182e-05\\
108	6.47881210634194e-05\\
109	6.35588055462044e-05\\
110	6.22610467669316e-05\\
111	6.08990075824378e-05\\
112	5.94764336369624e-05\\
113	5.79974207186029e-05\\
114	5.66359316290788e-05\\
115	5.52307862061174e-05\\
116	5.37715863587196e-05\\
117	5.2271348185029e-05\\
118	5.07191770254822e-05\\
119	4.91316931519048e-05\\
120	4.76425775202682e-05\\
121	4.58039666240695e-05\\
122	4.41105639316652e-05\\
123	4.24378321117747e-05\\
124	4.07365143150762e-05\\
125	3.90200091910375e-05\\
126	3.7280080520008e-05\\
127	3.55293727172801e-05\\
128	3.37632894336637e-05\\
129	3.20006935534689e-05\\
130	3.02399896602856e-05\\
131	2.84905394264401e-05\\
132	2.67608423071275e-05\\
133	2.73785103355425e-05\\
134	2.33459925275681e-05\\
135	2.16731417848615e-05\\
136	2.81461884772377e-05\\
137	1.84070251847968e-05\\
138	1.6816472404757e-05\\
139	1.52636172893184e-05\\
140	2.78674177702057e-05\\
141	1.226876436633e-05\\
142	2.93319468996096e-05\\
143	9.44704066813129e-06\\
144	8.1060706864089e-06\\
145	2.89880116408041e-05\\
146	1.03265949067494e-05\\
147	1.11191439790866e-05\\
148	4.44490746411214e-05\\
149	1.23064035240095e-06\\
150	1.43325681032432e-05\\
151	1.43222675043075e-05\\
152	1.55729873612959e-05\\
153	3.88914282481435e-06\\
154	1.96943184025107e-05\\
155	3.88528541591662e-06\\
156	1.08510610175028e-06\\
157	2.24816200636341e-05\\
158	1.85478429706096e-05\\
159	5.30931650201081e-06\\
160	1.02062124668345e-06\\
161	9.9900427769084e-07\\
162	4.19795980502548e-06\\
163	3.43957087497719e-06\\
164	1.05253036428443e-06\\
165	2.50770966744667e-05\\
166	2.06647163132675e-05\\
167	6.27029138089195e-06\\
168	8.31138093071449e-07\\
169	8.02778502395629e-07\\
170	2.61620554615265e-05\\
171	2.1183900557465e-05\\
172	6.58853675026915e-06\\
173	6.91349448356249e-07\\
174	2.73409751154493e-05\\
175	5.51966623694322e-06\\
176	2.84213885704242e-05\\
177	6.22753380848973e-06\\
178	2.83939029180078e-05\\
179	2.07560684120709e-05\\
180	2.01798429804209e-05\\
181	1.05081790587898e-05\\
182	1.45740073571925e-05\\
183	1.95068014413877e-06\\
184	1.03456603880546e-06\\
185	1.10612542071639e-05\\
186	8.54363713758372e-06\\
187	2.81370892824857e-06\\
188	2.74904062902426e-07\\
189	3.34338729375251e-07\\
190	2.27043107758895e-07\\
191	2.02217348088674e-07\\
192	1.78039991749479e-07\\
193	1.5686195503628e-07\\
194	2.78254009886317e-05\\
195	5.79278316726405e-06\\
196	1.55078034915156e-06\\
197	2.89785444678341e-05\\
198	4.44449183532192e-05\\
199	1.08505420204352e-07\\
};
\addlegendentry{Acc ADMM, restart every 50}

\end{axis}
\end{tikzpicture}%
\caption{Norm of primal residual versus number of iterations for the decentralized synthesis problem with standard and accelerated ADMM with various restart methods.}
\label{primal}
\end{figure}

%% file: ADHS_hybrid_dual.tex
%
%
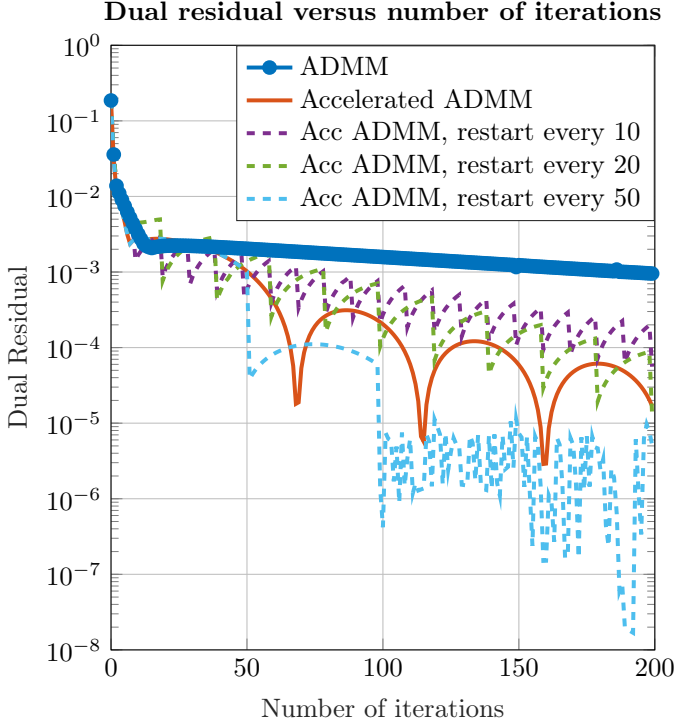
\begin{figure}[H]
\definecolor{mycolor1}{rgb}{0.00000,0.44700,0.74100}%
\definecolor{mycolor2}{rgb}{0.85000,0.32500,0.09800}%
\definecolor{mycolor3}{rgb}{0.92900,0.69400,0.12500}%
\definecolor{mycolor4}{rgb}{0.49400,0.18400,0.55600}%
\definecolor{mycolor5}{rgb}{0.46600,0.67400,0.18800}%
\definecolor{mycolor6}{rgb}{0.30100,0.74500,0.93300}%
\begin{tikzpicture}
\begin{axis}[%
width=2.817in,
height=3.156in,
at={(2.167in,0.898in)},
scale only axis,
xmin=0,
xmax=200,
yminorticks=true,
ymajorticks=true,
xmajorgrids=true,
ymajorgrids=true,
legend style={at={(1,1)},anchor=north east,font=\fontsize{10}{10}\selectfont},
xlabel style={font=\color{white!15!black}},
xlabel={Number of iterations},
ymode=log,
ymin=1e-08,
ymax=1,
yminorticks=true,
ylabel style={font=\color{white!15!black}},
ylabel={Dual Residual},
axis background/.style={fill=white},
title style={font=\bfseries},
title={Dual residual versus number of iterations},
legend style={legend cell align=left, align=left, draw=white!15!black}
]
\addplot [color=mycolor1, line width=1.5pt, draw=none, mark=*, mark options={solid, mycolor1}]
  table[row sep=crcr]{%
0	0.186075648451696\\
1	0.0358919761438512\\
2	0.0138582131543014\\
3	0.0110042925925947\\
4	0.00901825736693894\\
5	0.00748401727221372\\
6	0.0062406850288058\\
7	0.00523183778240891\\
8	0.0044066330884448\\
9	0.00373226815312278\\
10	0.00318162474134548\\
11	0.00273445776886638\\
12	0.00236274384549842\\
13	0.00220749462436057\\
14	0.00213138763162555\\
15	0.00207957875031881\\
16	0.00212724458823036\\
17	0.00216198008276528\\
18	0.00218766907091934\\
19	0.00220613328914662\\
20	0.00221909159910153\\
21	0.00222769653431179\\
22	0.00223310988491402\\
23	0.00223581953446487\\
24	0.00223647753677445\\
25	0.00223563112642046\\
26	0.00223270054332852\\
27	0.0022289081711083\\
28	0.00222417015047283\\
29	0.0022186289732766\\
30	0.00221240867326825\\
31	0.00220561229112903\\
32	0.00219831088885422\\
33	0.00219058975437891\\
34	0.00218252529159992\\
35	0.00217418026800516\\
36	0.00216556777053528\\
37	0.00215661795827307\\
38	0.0021474651922954\\
39	0.00213873696345304\\
40	0.0021292904971761\\
41	0.00211967689931734\\
42	0.00210998967122973\\
43	0.00210016499022443\\
44	0.00209038619351189\\
45	0.00208058083251952\\
46	0.00207073202217892\\
47	0.00206086714192556\\
48	0.00205099293124772\\
49	0.00204111383816852\\
50	0.00203123343670024\\
51	0.00202135638201485\\
52	0.00201147207163818\\
53	0.00200160390333547\\
54	0.00199172140205837\\
55	0.00198183441873784\\
56	0.00197203913780862\\
57	0.00196226463452812\\
58	0.00195251476456162\\
59	0.00194274315346059\\
60	0.00193305224325978\\
61	0.00192340678171662\\
62	0.00191378389198991\\
63	0.00190419959911828\\
64	0.00189465118760018\\
65	0.00188513965633697\\
66	0.00187566609823503\\
67	0.00186623110922354\\
68	0.00185683753702368\\
69	0.00184749126358574\\
70	0.00183817934737465\\
71	0.00182890885760562\\
72	0.00181967986116977\\
73	0.00181049925637192\\
74	0.00180135262276385\\
75	0.00179224804878714\\
76	0.0017831867438001\\
77	0.00177416852532594\\
78	0.0017651927996491\\
79	0.00175626022762912\\
80	0.00174737026513705\\
81	0.00173852339110223\\
82	0.00172971456680707\\
83	0.00172095622847801\\
84	0.00171224632483558\\
85	0.00170356854308859\\
86	0.00169495564618982\\
87	0.00168634449562186\\
88	0.00167779666957316\\
89	0.00166929800194254\\
90	0.0016608424389303\\
91	0.00165234866210878\\
92	0.00164396033129727\\
93	0.00163563687358917\\
94	0.00162735439149855\\
95	0.00161911257550318\\
96	0.00161090846361023\\
97	0.00160275032602847\\
98	0.00159462980785666\\
99	0.00158653922171001\\
100	0.00157850016243727\\
101	0.00157050175243384\\
102	0.00156261136418684\\
103	0.00155468470246145\\
104	0.00154678995835433\\
105	0.00153893697128196\\
106	0.00153111384200783\\
107	0.00152333908077759\\
108	0.001515619988046\\
109	0.00150793983487746\\
110	0.0015002891317267\\
111	0.00149269086893961\\
112	0.00148512301174647\\
113	0.00147759282273774\\
114	0.00147010109881487\\
115	0.00146264733464441\\
116	0.00145523081869406\\
117	0.00144785297450697\\
118	0.00144050970895027\\
119	0.00143320462708088\\
120	0.00142593712964975\\
121	0.00141870400645314\\
122	0.00141150897941061\\
123	0.00140435025592625\\
124	0.00139722846965273\\
125	0.00139014032321366\\
126	0.0013830882966348\\
127	0.00137607092021608\\
128	0.00136906487316985\\
129	0.00136214766485221\\
130	0.00135526613280165\\
131	0.00134837201831456\\
132	0.00134154114642246\\
133	0.00133473199478573\\
134	0.00132796553079695\\
135	0.00132123127004088\\
136	0.00131453036063781\\
137	0.00130786288322654\\
138	0.00130122783883896\\
139	0.0012946349891009\\
140	0.00128806167679765\\
141	0.00128152912405502\\
142	0.00127503089818501\\
143	0.0012685648954448\\
144	0.00126212914875973\\
145	0.00125572755769421\\
146	0.00124936687705335\\
147	0.00124302625954437\\
148	0.00123672547665105\\
149	0.00116308923028324\\
150	0.0012318879505015\\
151	0.00122280648069811\\
152	0.00121480642507243\\
153	0.00120744411128781\\
154	0.00120052131159371\\
155	0.00119390940174963\\
156	0.00118752582618065\\
157	0.00118129307191411\\
158	0.00117518764745406\\
159	0.00116918671048932\\
160	0.0011632009754356\\
161	0.00115731824727466\\
162	0.00115146123543258\\
163	0.00114565005333384\\
164	0.00113990019002236\\
165	0.00113412539623436\\
166	0.00112839501945039\\
167	0.00112270195319519\\
168	0.00111703543062552\\
169	0.00111140171130741\\
170	0.00110579222021484\\
171	0.00110020949531592\\
172	0.00109465498351452\\
173	0.00108912496494951\\
174	0.0010836284578743\\
175	0.0010781529166885\\
176	0.00107270437299173\\
177	0.00106728180069005\\
178	0.00106188568973198\\
179	0.00105651197147566\\
180	0.00105116946701522\\
181	0.00104585702696064\\
182	0.00104056213766152\\
183	0.00103530129615259\\
184	0.00103005708558862\\
185	0.00102484970211636\\
186	0.00107781691623443\\
187	0.00100739916093706\\
188	0.00100478130357518\\
189	0.00100129125118553\\
190	0.000997290443973176\\
191	0.000992968554934526\\
192	0.000988430355049252\\
193	0.000983739129865879\\
194	0.000978968085503629\\
195	0.000974126635529926\\
196	0.000969271748449629\\
197	0.000964407445170457\\
198	0.000959544356043285\\
199	0.000954691380232964\\
};
\addlegendentry{ADMM}

\addplot [color=mycolor2, line width=1.5pt]
  table[row sep=crcr]{%
0	0.186075648451696\\
1	0.0303821469144503\\
2	0.0118777268041626\\
3	0.00889881490269405\\
4	0.00653460464528318\\
5	0.00462982149983636\\
6	0.0031263521188918\\
7	0.00237176662121691\\
8	0.00259446529596231\\
9	0.00273314824073139\\
10	0.00276416656766667\\
11	0.0027512120955906\\
12	0.00273222495263368\\
13	0.00272400000760406\\
14	0.00272893112188236\\
15	0.00274070790669303\\
16	0.00275297552000909\\
17	0.00275950210996596\\
18	0.00275720961755044\\
19	0.00274542962811287\\
20	0.00272531917686928\\
21	0.00269898275102499\\
22	0.00266841885582367\\
23	0.00263533715991947\\
24	0.00260051163441412\\
25	0.00256405941104093\\
26	0.00252564169611491\\
27	0.00248479033845916\\
28	0.00242818877117179\\
29	0.00239567412584384\\
30	0.00234565519105233\\
31	0.00230246281558279\\
32	0.0022341840092959\\
33	0.00217593272793116\\
34	0.00211637937110081\\
35	0.00205537131815007\\
36	0.00199265253835181\\
37	0.00192805499008828\\
38	0.00186157693137207\\
39	0.00179348632920158\\
40	0.00172451089005235\\
41	0.00165460182647802\\
42	0.00158405749223102\\
43	0.00151303227938939\\
44	0.00144187666740977\\
45	0.00137075837139828\\
46	0.00129977181837863\\
47	0.00122903041238545\\
48	0.00115863421238691\\
49	0.00108870695047689\\
50	0.00101939046987855\\
51	0.000950832758734922\\
52	0.000883194562418097\\
53	0.000816639904278819\\
54	0.000751253750979286\\
55	0.000687231712638738\\
56	0.000624683367511989\\
57	0.000563708719396897\\
58	0.000504402022934431\\
59	0.000446881173118512\\
60	0.000391153139767342\\
61	0.000337324900043734\\
62	0.000285457386551654\\
63	0.00023559784423785\\
64	0.000187791699037915\\
65	0.000142087347208904\\
66	9.85290766340364e-05\\
67	5.71608440658493e-05\\
68	1.80162959358682e-05\\
69	1.88735986369376e-05\\
70	5.34842613621758e-05\\
71	8.58033378247733e-05\\
72	0.000115829860182731\\
73	0.000143573690217906\\
74	0.000169048217064124\\
75	0.000192282697219217\\
76	0.000213305965434332\\
77	0.000232159345469736\\
78	0.000248880626010141\\
79	0.000263517573011404\\
80	0.00027611564816738\\
81	0.000286709549334565\\
82	0.000295404001568721\\
83	0.000302224652050552\\
84	0.000307237780622047\\
85	0.000310505179006431\\
86	0.000312106483539849\\
87	0.000312111932183495\\
88	0.000310600338903654\\
89	0.000307667535408607\\
90	0.000303362773113915\\
91	0.000297787875380046\\
92	0.00029102873974109\\
93	0.000283161444257527\\
94	0.000274262969050129\\
95	0.000264415341888217\\
96	0.000253704478785188\\
97	0.000242218374052745\\
98	0.00023004688438837\\
99	0.000217279433631157\\
100	0.000204006813573979\\
101	0.000190307539077553\\
102	0.000176259737941303\\
103	0.000161946455174437\\
104	0.000147443559989148\\
105	0.000132823070301174\\
106	0.00011815420782505\\
107	0.000103504372341557\\
108	8.89389499134697e-05\\
109	7.45214169179929e-05\\
110	6.03124271283086e-05\\
111	4.63700799148134e-05\\
112	3.27489669635341e-05\\
113	1.94996637528095e-05\\
114	6.66946787951368e-06\\
115	5.69817301520601e-06\\
116	1.75641399604091e-05\\
117	2.88924937562946e-05\\
118	3.96511348459696e-05\\
119	4.98112148663652e-05\\
120	5.93474373631684e-05\\
121	6.82383655268883e-05\\
122	7.64653431172377e-05\\
123	8.4010698400854e-05\\
124	9.08702336391744e-05\\
125	9.70326107339185e-05\\
126	0.000102492265512848\\
127	0.000107248100814079\\
128	0.000111303718854239\\
129	0.000114662080293485\\
130	0.000117333390981878\\
131	0.000119329334968088\\
132	0.000120665103890694\\
133	0.000121354145633422\\
134	0.000121416409618475\\
135	0.000120873303646459\\
136	0.000119747125916036\\
137	0.000118062691390336\\
138	0.000115846481147504\\
139	0.000113126366493395\\
140	0.000109932932561176\\
141	0.000106294982170243\\
142	0.000102240463455751\\
143	9.78114739345241e-05\\
144	9.30372376799993e-05\\
145	8.79489022219304e-05\\
146	8.25823589552932e-05\\
147	7.69697800373311e-05\\
148	7.11534483484894e-05\\
149	6.51624922374451e-05\\
150	5.90290130367633e-05\\
151	5.27887352205564e-05\\
152	4.64765004775922e-05\\
153	4.012287374401e-05\\
154	3.37612156813636e-05\\
155	2.74239515427883e-05\\
156	2.11381698905195e-05\\
157	1.49342768352631e-05\\
158	8.83996475190705e-06\\
159	2.88257281238676e-06\\
160	2.91587778180141e-06\\
161	8.53050683386529e-06\\
162	1.39388035113991e-05\\
163	1.91198620331638e-05\\
164	2.40576768742284e-05\\
165	2.87328235630985e-05\\
166	3.31323080993635e-05\\
167	3.72409993609792e-05\\
168	4.1046937760969e-05\\
169	4.45401508740606e-05\\
170	4.77124299642711e-05\\
171	5.05573934965458e-05\\
172	5.30705003371599e-05\\
173	5.52485390852212e-05\\
174	5.70902328414485e-05\\
175	5.85945817453457e-05\\
176	5.97660546861174e-05\\
177	6.06053709954865e-05\\
178	6.11193829507122e-05\\
179	6.13120591046702e-05\\
180	6.11932309591001e-05\\
181	6.07702299786551e-05\\
182	6.00532579238593e-05\\
183	5.90536463487686e-05\\
184	5.77832032995359e-05\\
185	5.62576808734253e-05\\
186	5.44897466813878e-05\\
187	5.24943798002802e-05\\
188	5.02869259592116e-05\\
189	4.78837718843864e-05\\
190	4.53009538240287e-05\\
191	4.25606992093286e-05\\
192	3.9672472089935e-05\\
193	3.66600089767977e-05\\
194	3.35362676136877e-05\\
195	3.0322586813371e-05\\
196	2.70379904231398e-05\\
197	2.37025073258011e-05\\
198	2.03278911978006e-05\\
199	1.69343842067508e-05\\
};
\addlegendentry{Accelerated ADMM}

\addplot [color=mycolor4, dashed, line width=1.5pt]
  table[row sep=crcr]{%
0	0.186075648451696\\
1	0.0348999221819073\\
2	0.0147964517722459\\
3	0.0117420566606085\\
4	0.00901598785501945\\
5	0.00662513141730219\\
6	0.00461378655075478\\
7	0.00359492595967486\\
8	0.00402620948621611\\
9	0.00142101154376858\\
10	0.00168022290730856\\
11	0.00189572770988643\\
12	0.00208591187250671\\
13	0.00226323055637108\\
14	0.00243430797654383\\
15	0.00260288584771833\\
16	0.00276865519315268\\
17	0.00292988954176588\\
18	0.00308455877711434\\
19	0.000862717281220914\\
20	0.00108474833172059\\
21	0.00129175420385434\\
22	0.00148634308400529\\
23	0.00166979024063078\\
24	0.00184280531770422\\
25	0.00200580198472899\\
26	0.00215898164410084\\
27	0.00230268804616397\\
28	0.00243701216553661\\
29	0.000684104181674319\\
30	0.000863798284866418\\
31	0.00103261011424151\\
32	0.00119238983520226\\
33	0.00134400189845013\\
34	0.00148788160787474\\
35	0.00162427807369331\\
36	0.00175331316196116\\
37	0.00187505524457396\\
38	0.00198952887023151\\
39	0.000559708226722582\\
40	0.000708477585695417\\
41	0.000848927634278568\\
42	0.000982435125494563\\
43	0.00110964135283304\\
44	0.0012308785781386\\
45	0.00134624687901475\\
46	0.001455784791034\\
47	0.00155953620669015\\
48	0.00165742130687604\\
49	0.000467015042537153\\
50	0.000592025908557076\\
51	0.000710391690415859\\
52	0.000823241895888271\\
53	0.000931058045352888\\
54	0.00103405199470218\\
55	0.00113230461057526\\
56	0.00122582824044701\\
57	0.00131459888263902\\
58	0.00139858222069362\\
59	0.000394470435029552\\
60	0.000500560255140919\\
61	0.000601212074670142\\
62	0.00069735947216044\\
63	0.000789388626379643\\
64	0.000877460618576902\\
65	0.00096162808767358\\
66	0.00104188287198267\\
67	0.00109917691907183\\
68	0.00119391061967355\\
69	0.000336835640869007\\
70	0.000427342522251137\\
71	0.000513220972814954\\
72	0.000595358770021993\\
73	0.000674137310826131\\
74	0.000749709846936841\\
75	0.000822116389940431\\
76	0.000891321376352247\\
77	0.000957266905560575\\
78	0.00101987827973445\\
79	0.000288112334902325\\
80	0.000366048751816704\\
81	0.000440178575941913\\
82	0.000511166822535651\\
83	0.000579280400846174\\
84	0.000644622768695036\\
85	0.000722700128079748\\
86	0.000763928270414015\\
87	0.000821850328974803\\
88	0.000877207918067859\\
89	0.000248627416748291\\
90	0.000315858599098687\\
91	0.000379816105152254\\
92	0.000441096678272288\\
93	0.000499940398354193\\
94	0.000556447267737041\\
95	0.000610639848226389\\
96	0.000662501013551423\\
97	0.000711991186778197\\
98	0.000759060731396009\\
99	0.000214539181662533\\
100	0.000272780817489289\\
101	0.000328267450485852\\
102	0.000381486077414734\\
103	0.000432628726681549\\
104	0.000481765670348004\\
105	0.000528908185899305\\
106	0.000574037402413621\\
107	0.000617119248522927\\
108	0.000658111937487385\\
109	0.000186064282447741\\
110	0.000236637977024867\\
111	0.000284846241207133\\
112	0.000331110017308118\\
113	0.0003755935883254\\
114	0.000418357718819068\\
115	0.000459398493616894\\
116	0.000498709524008497\\
117	0.000536263020327584\\
118	0.000572017992891645\\
119	0.000161761080188165\\
120	0.000205775141445216\\
121	0.000247754059696417\\
122	0.000288054546504125\\
123	0.000326825294892728\\
124	0.000364109649213684\\
125	0.000399888571651736\\
126	0.000434217846335097\\
127	0.000467025521302423\\
128	0.000498276819526302\\
129	0.000140945167933111\\
130	0.000179328989915866\\
131	0.000215950017710544\\
132	0.000251124444418365\\
133	0.000284974773112187\\
134	0.000317544429197941\\
135	0.000348791446098159\\
136	0.000378777414211426\\
137	0.000407468390030094\\
138	0.000434807071101487\\
139	0.000123003662073948\\
140	0.000156527020013094\\
141	0.000188521618694072\\
142	0.000219263302223958\\
143	0.000248857870231582\\
144	0.000277342506023157\\
145	0.00030472036777062\\
146	0.000330977287406582\\
147	0.000356091491997463\\
148	0.000380031557397134\\
149	0.000107525587548438\\
150	0.000136851790989024\\
151	0.000164847275401089\\
152	0.000191755993772481\\
153	0.000217674628163433\\
154	0.000242622725219699\\
155	0.000266605981817262\\
156	0.000289612906522534\\
157	0.000311623948808229\\
158	0.000332620827316782\\
159	9.41260548258803e-05\\
160	0.000119810760705182\\
161	0.000144338797391257\\
162	0.000167919604652036\\
163	0.000190633090584792\\
164	0.00021250639169734\\
165	0.000233541601070263\\
166	0.000253726684100314\\
167	0.000273042810952655\\
168	0.000291468436719825\\
169	8.24867520615875e-05\\
170	0.000105007832056303\\
171	0.00012651971370874\\
172	0.000147205718978206\\
173	0.000167136459774812\\
174	0.000186331783748489\\
175	0.000204796166240334\\
176	0.000222518635759813\\
177	0.000239483106331367\\
178	0.000255669252002562\\
179	7.23623224114936e-05\\
180	9.21284693527534e-05\\
181	0.000111012270045681\\
182	0.00012917133218153\\
183	0.000146664756750366\\
184	0.000163526600433907\\
185	0.000179751762187053\\
186	0.000195328730153414\\
187	0.000210241556850184\\
188	0.00022447319721753\\
189	6.35395589291093e-05\\
190	8.0901972683875e-05\\
191	9.74924313144474e-05\\
192	0.000113452022630817\\
193	0.000128834673209601\\
194	0.000143657290318466\\
195	0.000157920278114776\\
196	0.000171616036442718\\
197	0.000184730604305895\\
198	0.000197249328449483\\
199	5.58368038574759e-05\\
};
\addlegendentry{Acc ADMM, restart every 10}

\addplot [color=mycolor5, dashed, line width=1.5pt]
  table[row sep=crcr]{%
0	0.186075648451696\\
1	0.0348999221819073\\
2	0.0147964517722459\\
3	0.0117420566606085\\
4	0.00901598785501945\\
5	0.00662513141730219\\
6	0.00461378655075478\\
7	0.00359492595967486\\
8	0.00402620948621611\\
9	0.00433174804072225\\
10	0.0044652188242365\\
11	0.00452230991551697\\
12	0.00456357427953557\\
13	0.00461777439805985\\
14	0.0046904102046112\\
15	0.00477184962400351\\
16	0.00485167983550668\\
17	0.00491909524398769\\
18	0.00496844520930244\\
19	0.000767324276311776\\
20	0.000960487061158481\\
21	0.00113961294490343\\
22	0.00130771465749302\\
23	0.00146641959016032\\
24	0.00161659682382382\\
25	0.00175863843660622\\
26	0.00189266205074011\\
27	0.00201862859826889\\
28	0.00213646176968663\\
29	0.00224601031019359\\
30	0.0023474189824404\\
31	0.0024407174324087\\
32	0.00252599530883324\\
33	0.00260332780157417\\
34	0.00267290732119028\\
35	0.00273483528215943\\
36	0.0027892445588823\\
37	0.00283626466690553\\
38	0.00287599167068241\\
39	0.000422708534750428\\
40	0.000535103868998134\\
41	0.000641219752983769\\
42	0.000742125079144536\\
43	0.00081688468752451\\
44	0.000933695241194566\\
45	0.00101965149319128\\
46	0.0011009711063376\\
47	0.00117813808902144\\
48	0.00125137741700002\\
49	0.00132070684205336\\
50	0.00138597498195243\\
51	0.00144697193333009\\
52	0.00150349950424131\\
53	0.00155541571255755\\
54	0.00162059859084438\\
55	0.00164154745548262\\
56	0.00168028076604438\\
57	0.0017150156300944\\
58	0.00174519940666652\\
59	0.000258071607727425\\
60	0.000327157314921638\\
61	0.000392622537041761\\
62	0.000455134172411349\\
63	0.000514984990987345\\
64	0.000572303230981256\\
65	0.000627130630601618\\
66	0.000679459456780097\\
67	0.000729254541422315\\
68	0.00077646913339811\\
69	0.000821052976066321\\
70	0.000862960335924875\\
71	0.000902150636229747\\
72	0.000938591183508836\\
73	0.00097224681055512\\
74	0.0010031049304709\\
75	0.00103115560473755\\
76	0.001056383780378\\
77	0.00107872256996706\\
78	0.00109837037354698\\
79	0.000162198225178494\\
80	0.000206072566937779\\
81	0.000247802880234083\\
82	0.000287764814199005\\
83	0.000326109280403011\\
84	0.000362888293994478\\
85	0.000398064574425027\\
86	0.000431760244939348\\
87	0.000463883238405488\\
88	0.000494392982239324\\
89	0.000523262615619246\\
90	0.000550468776279228\\
91	0.000575984110384507\\
92	0.000599784274286777\\
93	0.000621846282367256\\
94	0.000642143559897743\\
95	0.000660674345127838\\
96	0.000677415241091246\\
97	0.000692373459290396\\
98	0.000705524990098913\\
99	0.00010424273819588\\
100	0.000132542359972077\\
101	0.000159503479016996\\
102	0.000185361890030266\\
103	0.000210210930307072\\
104	0.000234084875977089\\
105	0.000256989698305567\\
106	0.000278916823018533\\
107	0.000299848915174762\\
108	0.000319766528891016\\
109	0.000338647923478116\\
110	0.000356470979387662\\
111	0.000373215439922324\\
112	0.000388860520977603\\
113	0.000403381089813325\\
114	0.000416764619150714\\
115	0.00042901851309737\\
116	0.000440121475543397\\
117	0.000450063274149796\\
118	0.000458839060224525\\
119	6.78336605978991e-05\\
120	8.62896665129589e-05\\
121	0.000103890415259534\\
122	0.000120788797760623\\
123	0.00013704402803685\\
124	0.000152677467979901\\
125	0.000167691337169901\\
126	0.00018207944060832\\
127	0.000195829990525362\\
128	0.000208928111124535\\
129	0.000221358150614148\\
130	0.000233104551109233\\
131	0.000244153304479164\\
132	0.000254489566026558\\
133	0.000264101722396818\\
134	0.000272978136903285\\
135	0.000281109562003076\\
136	0.000288488130707355\\
137	0.000295109310658851\\
138	0.000300968752880374\\
139	4.45075956282005e-05\\
140	5.66378047046408e-05\\
141	6.82145916573525e-05\\
142	7.93383845085283e-05\\
143	9.00462785904666e-05\\
144	0.000100352685326602\\
145	0.000110257988540862\\
146	0.000119758758715352\\
147	0.000128842327883518\\
148	0.000137501848424549\\
149	0.000145727031721976\\
150	0.000153507009265962\\
151	0.000160831239375224\\
152	0.000167689875211778\\
153	0.000174074173141357\\
154	0.000179975936702225\\
155	0.000185389140298221\\
156	0.000190311345951809\\
157	0.000194732458004809\\
158	0.000198649860888941\\
159	2.93823629943674e-05\\
160	3.74008261099332e-05\\
161	4.50581572770643e-05\\
162	5.24194294927325e-05\\
163	5.95104336372934e-05\\
164	6.63388446005734e-05\\
165	7.29049783414459e-05\\
166	7.92056060101945e-05\\
167	8.52432048736788e-05\\
168	9.09960333139876e-05\\
169	9.64600928891687e-05\\
170	0.000101629738135074\\
171	0.000106498932612508\\
172	0.000111063053073422\\
173	0.000115315896714906\\
174	0.000119254280144477\\
175	0.00012286145437183\\
176	0.000126146529399278\\
177	0.000129109800207928\\
178	0.000131741635140832\\
179	1.9489499758149e-05\\
180	2.48189867236647e-05\\
181	2.99045074699841e-05\\
182	3.47953954071194e-05\\
183	3.95084459271675e-05\\
184	4.40492663990457e-05\\
185	4.84179875441877e-05\\
186	5.26119670759107e-05\\
187	5.66273317212002e-05\\
188	6.04532609420642e-05\\
189	6.4101102767234e-05\\
190	6.75552924983135e-05\\
191	7.08113148875553e-05\\
192	7.38589849699599e-05\\
193	7.67017004851254e-05\\
194	7.93299807905313e-05\\
195	8.17452119428362e-05\\
196	8.39460236125032e-05\\
197	8.59294839968142e-05\\
198	8.769324235204e-05\\
199	1.29701182701878e-05\\
};
\addlegendentry{Acc ADMM, restart every 20}
\addlegendentry{Acc ADMM, restart every 50}
\addplot [color=mycolor6, dashed, line width=1.5pt]
  table[row sep=crcr]{%
0	0.186075648451696\\
1	0.0303821469144503\\
2	0.0118777268041626\\
3	0.00889881490269405\\
4	0.00653460464528318\\
5	0.00462982149983636\\
6	0.0031263521188918\\
7	0.00237176662121691\\
8	0.00259446529596231\\
9	0.00273314824073139\\
10	0.00276416656766667\\
11	0.0027512120955906\\
12	0.00273222495263368\\
13	0.00272400000760406\\
14	0.00272893112188236\\
15	0.00274070790669303\\
16	0.00275297552000909\\
17	0.00275950210996596\\
18	0.00275720961755044\\
19	0.00274542962811287\\
20	0.00272531917686928\\
21	0.00269898275102499\\
22	0.00266841885582367\\
23	0.00263533715991947\\
24	0.00260051163441412\\
25	0.00256405941104093\\
26	0.00252564169611491\\
27	0.00248479033845916\\
28	0.00242818877117179\\
29	0.00239567412584384\\
30	0.00234565519105233\\
31	0.00230246281558279\\
32	0.0022341840092959\\
33	0.00217593272793116\\
34	0.00211637937110081\\
35	0.00205537131815007\\
36	0.00199265253835181\\
37	0.00192805499008828\\
38	0.00186157693137207\\
39	0.00179348632920158\\
40	0.00172451089005235\\
41	0.00165460182647802\\
42	0.00158405749223102\\
43	0.00151303227938939\\
44	0.00144187666740977\\
45	0.00137075837139828\\
46	0.00129977181837863\\
47	0.00122903041238545\\
48	0.00115863421238691\\
49	0.00108870695047689\\
50	0.00101939046987855\\
51	3.82804945286469e-05\\
52	4.42792802802065e-05\\
53	5.00023186169313e-05\\
54	5.54702173620235e-05\\
55	6.06926254904111e-05\\
56	6.56776986261951e-05\\
57	7.0415803059029e-05\\
58	7.49109628617167e-05\\
59	7.91593539390125e-05\\
60	8.31571466641603e-05\\
61	8.68969825692009e-05\\
62	9.03734958107211e-05\\
63	9.35814793826708e-05\\
64	9.65080603356635e-05\\
65	9.91741260466639e-05\\
66	0.000101563605315661\\
67	0.000103673330898794\\
68	0.000105504446475319\\
69	0.000107062178053476\\
70	0.000108341746574685\\
71	0.00010935233146779\\
72	0.000110099686101951\\
73	0.000110588499716636\\
74	0.000110823750559256\\
75	0.000110810294958031\\
76	0.000110554202044276\\
77	0.000110061778248949\\
78	0.000109342656743908\\
79	0.000108395399067824\\
80	0.000107232822911621\\
81	0.000105855495188638\\
82	0.000104293339995458\\
83	0.000102541923882591\\
84	0.000100607207226986\\
85	9.85009085863285e-05\\
86	9.62277986683001e-05\\
87	9.37948125720483e-05\\
88	9.12230179329486e-05\\
89	8.85151949675967e-05\\
90	8.56873399226792e-05\\
91	8.27511215689311e-05\\
92	7.97166574906007e-05\\
93	7.65933208020311e-05\\
94	7.33911629802804e-05\\
95	7.01224071081462e-05\\
96	6.67966616650752e-05\\
97	6.34281370732499e-05\\
98	6.0019315013106e-05\\
99	1.34138672794091e-06\\
100	4.20202906712173e-07\\
101	7.40494394011416e-06\\
102	6.19093126152923e-06\\
103	5.89859926666333e-06\\
104	1.70173644198248e-06\\
105	4.890973199435e-06\\
106	8.85802039755382e-07\\
107	7.49346957953901e-06\\
108	1.48381202542556e-06\\
109	3.43375537473713e-06\\
110	6.16550301151261e-06\\
111	1.18578486517139e-06\\
112	1.60833253877626e-06\\
113	1.28431231879547e-06\\
114	1.32823252725908e-06\\
115	1.36852243935821e-06\\
116	7.09221067937367e-06\\
117	5.71309239492133e-06\\
118	1.46302993357698e-06\\
119	1.11931156039689e-05\\
120	8.9872334993138e-06\\
121	7.09342585190636e-06\\
122	1.54368254445914e-06\\
123	7.38021828748674e-06\\
124	2.41432672288989e-06\\
125	4.14167696562635e-06\\
126	6.36979104623768e-06\\
127	5.22075591899459e-06\\
128	5.21999689553287e-06\\
129	1.53219984406316e-06\\
130	1.51322688270765e-06\\
131	1.49720843006934e-06\\
132	6.84979182433229e-06\\
133	5.47959298310844e-06\\
134	1.74557522176029e-06\\
135	6.95147280429664e-06\\
136	1.61771808373903e-06\\
137	1.34073273265337e-06\\
138	1.30408578517435e-06\\
139	6.77575486311596e-06\\
140	1.35715644452076e-06\\
141	7.08018823082673e-06\\
142	1.6128321433713e-06\\
143	1.1026973889714e-06\\
144	6.92530544414968e-06\\
145	2.45891436650937e-06\\
146	2.63907030825251e-06\\
147	1.05162745027279e-05\\
148	2.59905249077725e-06\\
149	1.13142892677054e-05\\
150	7.32463786262131e-06\\
151	6.4963909341798e-06\\
152	1.43608516871603e-06\\
153	6.69717046329645e-06\\
154	1.24439595171034e-06\\
155	2.38587746592647e-07\\
156	6.63231516152942e-06\\
157	5.31159841286207e-06\\
158	1.4831514672926e-06\\
159	1.50600422279294e-07\\
160	1.51189755558358e-07\\
161	1.10946346825942e-06\\
162	8.96230504857375e-07\\
163	2.70794939340049e-07\\
164	6.37841132114193e-06\\
165	5.20171194009898e-06\\
166	1.56339808498537e-06\\
167	1.63827375290829e-07\\
168	1.75884573908768e-07\\
169	6.36568518677532e-06\\
170	5.11819419871556e-06\\
171	1.58134620099741e-06\\
172	1.63755514775884e-07\\
173	6.48323324793131e-06\\
174	1.30157372406989e-06\\
175	6.66652819694401e-06\\
176	1.45337411737637e-06\\
177	6.59464673272576e-06\\
178	4.79850488057776e-06\\
179	4.64466597467584e-06\\
180	2.40832207261971e-06\\
181	3.32647991143212e-06\\
182	4.43479742858342e-07\\
183	2.34307332102318e-07\\
184	2.49589187230819e-06\\
185	1.92091355258984e-06\\
186	6.30426150529602e-07\\
187	5.96515066064028e-08\\
188	7.44129639877552e-08\\
189	1.94704162937543e-08\\
190	1.85479932966779e-08\\
191	1.7691662633523e-08\\
192	1.68494328713956e-08\\
193	6.09948894849811e-06\\
194	1.26621599607316e-06\\
195	3.38038050245361e-07\\
196	6.29957871956535e-06\\
197	9.63607696233859e-06\\
198	7.68520770778421e-06\\
199	5.13028503508075e-06\\
};

\end{axis}

\end{tikzpicture}%
\caption{Norm of dual residual versus number of iterations for the decentralized synthesis problem with standard and accelerated ADMM with various restart methods.}
\label{dual}
\end{figure}